\documentclass[opre]{informs3nothing} % current default for manuscript submission to arxiv (take out INFORMS stuff)
\RequirePackage[OT1]{fontenc}
\RequirePackage{
amsmath,
times,
graphicx,
amsfonts,
amscd,
amssymb,
algorithm,
algorithmic
}
\OneAndAHalfSpacedXI
\usepackage[sort&compress,comma,authoryear]{natbib}
 \bibpunct[, ]{(}{)}{,}{a}{}{,}%
 \def\bibfont{\small}%
\TheoremsNumberedThrough     % Preferred (Theorem 1, Lemma 1, Theorem 2)
\ECRepeatTheorems

\EquationsNumberedThrough    % Default: (1), (2), ...
\newcommand{\paranth}[1]{\left(#1\right)}
\newcommand{\bracket}[1]{\left[#1\right]}
\newcommand{\curly}[1]{\left\{#1\right\}}

\numberwithin{equation}{section}
\theoremstyle{plain}
\newtheorem{thm}{Theorem}[section]
\newtheorem{dfn}{Definition}
\newtheorem{lem}{Lemma}
\newtheorem{cor}{Corollary}
\newtheorem{prop}{Proposition}
\newcommand{\G}{\mathcal{G}}
\newcommand{\V}{\mathcal{V}}
\newcommand{\E}{\mathcal{E}}

\newcommand{\D}{\mathcal{D}}
\newcommand{\Eb}{\mathbf{E}}
\newcommand{\N}{\mathcal{N}}
\newcommand{\bepsilon}{\varepsilon}
\newcommand{\Tin}{T}

\begin{document}

% Outcomment only when entries are known. Otherwise leave as is and
%   default values will be used.
%\setcounter{page}{1}
%\VOLUME{00}%
%\NO{0}%
%\MONTH{Xxxxx}% (month or a similar seasonal id)
%\YEAR{0000}% e.g., 2005
%\FIRSTPAGE{000}%
%\LASTPAGE{000}%
%\SHORTYEAR{00}% shortened year (two-digit)
%\ISSUE{0000} %
%\LONGFIRSTPAGE{0001} %
%\DOI{10.1287/xxxx.0000.0000}%

% Author's names for the running heads
% Sample depending on the number of authors;
% \RUNAUTHOR{Jones}
% \RUNAUTHOR{Jones and Wilson}
% \RUNAUTHOR{Jones, Miller, and Wilson}
% \RUNAUTHOR{Jones et al.} % for four or more authors
% Enter authors following the given pattern:
\RUNAUTHOR{Shah and Zaman}

% Title or shortened title suitable for running heads. Sample:
% \RUNTITLE{Bundling Information Goods of Decreasing Value}
% Enter the (shortened) title:
\RUNTITLE{Finding Rumor Sources}

% Full title. Sample:
% \TITLE{Bundling Information Goods of Decreasing Value}
% Enter the full title:
\TITLE{Finding Rumor Sources on Random Trees}

% Block of authors and their affiliations starts here:
% NOTE: Authors with same affiliation, if the order of authors allows,
%   should be entered in ONE field, separated by a comma.
%   \EMAIL field can be repeated if more than one author
\ARTICLEAUTHORS{%
\AUTHOR{Devavrat Shah}
\AFF{Department of Electrical Engineering and Computer Science, Massachusetts
Institute of Technology, Cambridge, MA 02139, \EMAIL{devavrat@mit.edu}} %, \URL{}}

\AUTHOR{Tauhid Zaman}
\AFF{Sloan School of Management, Massachusetts
Institute of Technology, Cambridge, MA 02139, \EMAIL{zlisto@mit.edu}}
% Enter all authors
} % end of the block

\ABSTRACT{
We consider the problem of detecting the source of a rumor which has spread 
in a network using only observations about which set of nodes are infected with
 the rumor and with no information as to \emph{when} these nodes became infected. 
In a recent work \citep{ref:rc} this rumor source detection
problem was introduced and studied. 
The authors proposed the graph score function {\em rumor centrality} as an estimator for detecting the source. 
They establish it to be the maximum likelihood estimator with respect to 
the popular Susceptible Infected (SI) model with exponential spreading times
for regular trees. They showed that as the size of the infected graph increases,
for a path  graph (2-regular tree), the probability of source detection goes to $0$ 
while for $d$-regular trees with $d \geq 3$ the probability of detection, 
say $\alpha_d$, remains bounded away from $0$ and is less than $1/2$. However, their 
results stop short of providing insights for the performance of the
rumor centrality estimator in more general settings
such as irregular trees or the SI model with non-exponential spreading times. 

This paper overcomes this limitation and establishes the effectiveness of rumor centrality 
for source detection for generic random trees and the SI model with a generic spreading 
time distribution. The key result is an interesting connection between a 
continuous time branching process and the effectiveness of rumor 
centrality. Through this, it is possible to quantify the detection probability 
precisely. As a consequence, we recover all previous results  as a special 
case and obtain a variety of novel results including the {\em universality} 
of rumor centrality in the context of tree-like graphs  and the SI model with a 
generic spreading time distribution.
}

 \HISTORY{}

\maketitle

%%%%%%%%%%%%%%%%%%%%%%%%%%%%%%%%%%%%%%%%%%%%%%%%%%%%%%%%%%%%%%%%%%%%%%
\section{Introduction}
Imagine someone starts a rumor which then spreads through a social network.
After the rumor has spread for a long amount of time, we observe this 
network of rumor infected individuals.  
We only know who has heard the rumor and the underlying network structure.  
No information is given about \emph{when} the people heard the rumor.  
Our goal is to use only this information to discover the source
of the rumor.   
 
This \emph{rumor source detection problem} is very general and arises in many 
different contexts.  For example, the rumor could be a computer virus spreading through
 the Internet, a contagious disease infecting a human population, or a trend or new product
diffusing through a social network.  In each of these different scenarios, detection of the source is 
of great interest. 
One would naturally like to find the originator
of a malicious cyber-attack.  Detecting the source of a viral
epidemic would aid with the development of effective vaccination, quarantine and
prevention strategies.  
In social networks, sources of rumors, trends or new product adoption
may be effective at disseminating information,
and their identification would be of interest to companies wishing
to develop viral marketing campaigns. 

Detection of the source is made challenging in each of these situations by
the fact that one may not
have information regarding the time of the infection or adoption.  
For example, if the computer virus remains dormant and then upon activation 
renders the system inoperable, it may
not be possible to determine when the machine was infected.     For contagious
diseases, determining exactly when a person became infected can be difficult due to 
lack of sufficient data.  Rather, only a broad time window of when the infection occurred 
may be known.  For trends or new product adoption one may be able to determine the exact
time of adoption if this occurs through a social network such as Facebook which
records the time of each user's activity.  However, there can be situations 
where people do not share the fact that they have adopted until much after they have done so, 
making it difficult to pinpoint precisely when the adoption occurred.

Given the wide ranging applications, it begs to understand the fundamental limitations of the source
detection problem. Concretely, there are two key questions
that need to be addressed. First, how does one actually  construct the rumor source estimator?  
Since no information about infection times is given, a rumor source estimator would need to extract 
all information about the identity of the source using only the structure of the rumor infected network, 
but it is not obvious in what manner.  Second, what are the fundamental limits to this rumor source 
detection problem?  In particular, how accurately can one find the rumor source, what is the magnitude 
of  errors made in this detection, and how does the network structure affect one's ability to find the rumor source?  

\subsection{Related Work}
Rumor spreading was originally studied in the context of epidemiology in order
to predict, control, and prevent the spread of infectious diseases.  The epidemiological
models for the spread of disease
generally consisted of individuals that could be in one of three states
: susceptible, infected, or recovered.  In the susceptible-infected-recovered or 
SIR model all three states are
allowed, but there are variants such as the SI model which only consider
susceptible and infected individuals.
Daniel Bernoulli developed the first differential equation models for the spread of
a disease \citep{ref:epdiffeq1}.  Modern differential equation models
were introduced in \citep{ref:epdiffeq2} and later expanded in \citep{ref:epdiffeq3} and \citep{ref:epdiffeq4}.
These models provided insight into the disease spreading dynamics, but they were
very {coarse} and made several simplifying assumptions about human populations.  
The next level of modeling involved taking into account the network over which the disease spread.
Contact network modeling was able to capture in greater detail the specific manner
by which disease spread.  These models have allowed researchers to understand how
 the network structure affects the ability of a disease to become an epidemic 
(\cite{ref:ep1}, \cite{ref:ep2}, \cite{ref:ep3}).  The insights obtained from modeling
disease spreading at a network level have allowed epidemiologists to
develop vaccination and quarantine strategies
to control modern viral epidemics (\cite{ref:epcontrol1}, \cite{ref:epcontrol2}, \cite{ref:epcontrol3},  
\cite{ref:epcontrol4}, \cite{ref:epcontrol5}, \cite{ref:epcontrol6}, \cite{ref:epcontrol7}, \cite{ref:epcontrol8}).  

The network models developed for disease propagation have
found application in the context of online social networks.  In  
\cite{ref:networkad1}, \cite{ref:networkad2}, and \cite{ref:networkad3}
optimization methods were applied to network models to select
the best set of users to seed with a new product or information in order to 
maximize its spread in a social network.
This work is complementary to that in epidemiology, where {the goal}
is to prevent the spread of a viral outbreak, not accelerate it.
Another interesting line of work has focused on 
using the spread of rumors in a social network to reconstruct
the unknown network structure (\cite{ref:rumorRecon1}, \cite{ref:rumorRecon2}, \cite{netrapalli2012learning}).

Controlling the spread of a rumor, whether it be a contagious
disease or the adoption of a new product, has been the main focus of a large amount of research,
but the question of identifying
the source of the rumor has been largely overlooked.
A problem located at the intersection of probability theory and
information theory recently emerged which is thematically related 
to rumor source detection.  It is known as the \emph{reconstruction problem} and the
goal is to estimate the information possessed by a source based on
noisy observations about this information as it propagates through a network. 
There are interesting similarities between the two problems: the 
signal of interest, the information of the source (for the reconstruction problem) and the 
rumor source itself (for the rumor source detection problem) are extremely 
`low-dimensional'.  However, the observations for each problem, 
the noisy versions of the information (reconstruction) and infected nodes 
(rumor source detection), lie in a very `high-dimensional' setting.  
This makes estimation and
detection quite challenging.
%\footnote{The phrase {\em finding a needle in a 
%haystack} seems quite appropriate.}.  
It is not surprising that 
results for the reconstruction problem, even for tree or tree-like 
graphs, have required sophisticated mathematical techniques  
(\cite{ref:recon1}, \cite{ref:recon2}, \cite{ref:recon3}).  
Therefore, one would expect similar types of challenges
for the rumor source detection problem, which involves not
estimating information at a known source, but rather finding 
the source itself among a large number of vertices in a network.

The rumor source detection problem was first formally posed 
and studied in \cite{ref:rc}. The authors proposed 
a graph-score function called \emph{rumor centrality}  
as an estimator for the rumor source.  
They showed that the node with maximal rumor centrality 
is the maximum likelihood (ML) estimate of the source
for rumor spreading on regular trees
under the SI model with homogeneous exponential spreading times. 
They demonstrated the effectiveness of this estimator by establishing that  
the rumor source is found with strictly positive probability for regular 
trees and geometric trees under this setting.  The model and precise 
results from \cite{ref:rc} are described in Section \ref{sec:RCE}. 

While this work laid the foundations of the rumor source detection problem, the 
results had some key limitations.  First, they do not quantify the exact detection 
probability, say $\alpha_d$, for $d$-regular graphs, for the proposed 
ML estimator other than $\alpha_2 = 0$, $\alpha_3 = 0.25$ 
and $0 < \alpha_d \leq 0.5$ for $d \geq 4$ for the SI model with 
exponential spreading times.  Second, the results do not quantify 
the magnitude of the error in the event of not being able to 
identify the source.  Third, the results do not provide any insights 
into how the estimator behaves for rumor spreading on generic heterogeneous tree 
(or tree-like) graphs under the SI model with 
a generic spreading time distribution.  

%%%%%%%%%%%%%%%%%%%%%%%%%%%%%%%%%%%%%%%%%%%%%%%%%%%%%%%%%%
%%%%%%%%%%%%%%%%%%%%%%%%%%%%%%%%%%%%%%%%%%%%%%%%%%%%%%%%%%
%%%%%%%%%%%%%%%%%%%%%%%%%%%%%%%%%%%%%%%%%%%%%%%%%%%%%%%%%%
%%%%%%%%%%%%%%%%%%%%%%%%%%%%%%%%%%%%%%%%%%%%%%%%%%%%%%%%%%

\subsection{Summary of Results}

The primary reason behind the limitations of the results in \cite{ref:rc} is the fact that the analytic method 
employed there is quite specific to regular trees with homogeneous exponential spreading times. To overcome these
limitations, as the main contribution of this work we introduce a novel analysis method that utilizes connections 
to the classical Markov branching process (MBP) (equivalently, a generalized Polya's urn (GPU)). As a 
consequence of this, we are able to quantify the probability of the error event precisely and thus eliminate the 
shortcomings of the prior work. 

Our results in this work collectively establish that, even though, rumor centrality is an ML estimator only for regular trees and the SI model with exponential spreading times, it is universally effective with respect to heterogeneity in the tree structure and spreading time distributions.  It's effectiveness for generic random trees immediately implies its utility for finding sources in sparse random graphs that are locally tree-like. Examples include Erdos-Renyi and random regular graphs. A brief discussion to this effect can be found in Section \ref{ssec:randgraph}. 

The following is a summary of our main results (see Section \ref{sec:results} for 
precise statements): 

\begin{itemize}
\item[1.] {\em Regular trees, SI model with exponential spreading times}: 

We characterize $\alpha_d$, the detection 
probability for $d$-regular trees, for all $d$. Specifically, for $d\geq 3$ 
$$ \alpha_d = d I_{1/2}\Big(\frac{1}{d-2}, \frac{d-1}{d-2}\Big) - (d-1). $$ 
In above $I_{x}(a, b)$ is the incomplete beta function with parameters $a, b$ evaluated at $x \in [0,1]$ (see \eqref{eq:beta}). 
This implies that $\alpha_d > 0$ for $d \geq 3$, $\alpha_3 = 0.25$, and $\alpha_d \to 1-\ln 2$ as $d\to\infty$.
Further, we show that the probability of rumor centrality estimating the $k^{th}$ infected node as the source
decays as $\exp(-\Theta(k))$. The precise results are stated as Theorem \ref{thm:regTree}, Corollaries \ref{corollary:randRegInf} and 
\ref{corollary:randReg3}.
 
% 
% is related to the probability of a certain event involving a Beta distribution with parameters dependent on $d$. We establish
% that $\alpha_d$ is non-decreasing in $d$ and $\alpha_d \uparrow 1-\ln 2$ as $d\to \infty$. 

\item[2.] {\em Generic random trees, SI model with exponential spreading times}: 
For generic random trees (see Section \ref{ssec:randTree} 
for precise definition) which are expanding, we establish that there is strictly positive probability of correct detection using
rumor centrality. Furthermore, the probability of rumor centrality estimating the $k^{th}$ infected node as the source
decays as $\exp\paranth{-\Theta(k)}$. The precise results are stated as Theorem \ref{thm:randtree} and Theorem \ref{thm:randTreeError}.

\item[3.] {\em Geometric trees, SI model with generic spreading times}: 

For any geometric tree (see 
Section \ref{ssec:geom} for precise definition), we establish that the probability of correct detection goes to $1$ as the
number of infected nodes increases. The precise result is stated as Theorem \ref{thm:geom}.

\item[4.] {\em Generic random trees, SI model with generic spreading times}: 

For generic expanding random trees with generic spreading times
(see Section \ref{ssec:randTree} for definition), we establish that the probability of correct source detection remains bounded 
away from $0$. The precise result is stated as Theorem \ref{thm:randtree}.

\end{itemize}

%%%%%%%%%%%%%%%%%%%%%%%%%%%%%%%%%%%%%%%%%%%%%%%%%%%%%%%%%%%%%%%%%%%%%%%
%%%%%%%%%%%%%%%%%%%%%%%%%%%%%%%%%%%%%%%%%%%%%%%%%%%%%%%%%%%%%%%%%%%%%%%
%%%%%%%%%%%%%%%%%%%%%%%%%%%%%%%%%%%%%%%%%%%%%%%%%%%%%%%%%%%%%%%%%%%%%%%
\section{Model, Problem Statement and Rumor Centrality}\label{sec:RCE}

We start by describing the model and problem statement followed by a quick recall of the 
precise results from \cite{ref:rc}. In the process, we shall recall the definition of
rumor centrality and source estimation as introduced in \cite{ref:rc}.

%%%%%%%%%%%%%%%%%%%%%%%%%%%%%%%%%%%%%%%%%%%%%%%%%%%%%%%%%%%%%%%%%%%%%%%%%%%%%%%%%
\subsection{Model} \label{ssec:si}

Let $\G = (\V,\E)$ be a possibly infinite connected graph. Let $v \in \V$ be a rumor source 
from which a rumor starts spreading at time $0$. As per the classical Susceptible Infected 
(SI) model the rumor spreads in the graph. Specifically, each edge $e = (u_1,u_2)$ has a 
spreading time $S_e$ associated with it. If node $u_1$ gets infected at time $t_1$, then 
at time $t_1+S_e$ the infection spreads from $u_1$ to $u_2$. A node, once becoming infected,
remains infected. The spreading times associated with edges are independent random variables
\newcommand{\Real}{{\mathbb R}}
with identical distribution. Let $F: \Real \to [0,1]$ denote the cumulative density function of the
spreading time distribution. We shall assume that the distribution is non-negative valued, i.e.
$F(0) = 0$ and it is non-atomic at $0$, i.e. $F(0^+) = 0$. Since it is a cumulative density function, it
is non-decreasing and $\lim_{x\to\infty} F(x) = 1$. The simplest, homogeneous SI model has exponential
spreading times with parameter $\lambda > 0$ with $F(x) = 1-\exp(-\lambda x)$ for $x \geq 0$. In \cite{ref:rc},
the results were restricted to this homogeneous exponential spreading time setting. In this paper,
we shall develop results for arbitrary spreading time distributions consistent with the above assumptions.

Given the above spreading model, we observe the rumor infected graph $G(t) = (V(t),E(t))$ 
at some time $t > 0$.  To simplify our notation, we will refer to the time
dependent rumor infected graph at time $t$ simply as $G=(V,E)$. We do not know the value of $t$ or the realization of the spreading times on edges $e \in E$;
we only know the rumor infected nodes $V \subset \V$ and edges between them $E = V \times V \cap \E$. The 
goal is to find the rumor source (among $V$) given $G$. 
 
We note here that in this setting we do not observe the underlying graph $\mathcal G$.  This means
we do not observe edges on the boundary between infected and non-infected nodes.  
However, these boundary edges do provide additional information.
For example, if an infected node has a large number of 
uninfected neighbors, then it is likely that this node has 
not been infected for very long, otherwise more of its neighbors would be infected.
Intuitively, this would mean that it is less likely that this node is the source.
Our rumor source estimator, which we present next, does not require any knowledge
of $\mathcal G$, though our analysis of the estimator's performance will require
knowledge of the structure of $\mathcal G$.  We will find that
without observing $\mathcal G$, our rumor source estimator is still
able to perform well on a variety of graphs under general spreading models.

%%%%%%%%%%%%%%%%%%%%%%%%%%%%%%%%%%%%%%%%%%%%%%%%%%%%%%%%%%%%%%%%%%%%%%%%%

%%%%%%%%%%%%%%%%%%%%%%%%%%%%%%%%%%%%%%%%%%%%%%%%%%%%%%%%%%%%%%%%%%%%%%

\subsection{Rumor Centrality: An Estimator}

\begin{figure}\centering
	\includegraphics[scale=1]{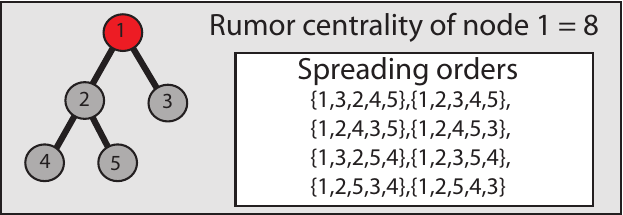}
	\caption{Example of rumor centrality calculation for a $5$ node network.  The rumor centrality of node $1$ is $8$ because there are $8$ spreading orders that it can originate, which are shown in the figure.}
	\label{fig:rcexample}
\end{figure}

To solve the rumor source detection problem, the notion of rumor centrality was introduced in \cite{ref:rc}.
Rumor centrality is a `graph score' function. That is, it takes $G = (V,E)$ as input
and assigns a non-negative number or score to each of the vertices. Then the estimated source
is the one with maximal (ties broken uniformly at random) score or rumor centrality. The 
node with maximal rumor centrality is called the `rumor center' (which is
also the estimated source) with ties broken uniformly at random. We start with the precise description 
of rumor centrality for a tree\footnote{ We shall call an undirected graph a {\em tree} 
if it is connected and it does not have any cycles.} graph $G$: the rumor
centrality of node $u\in V$ with respect to $G=(V,E)$ is
%%%
\begin{align}\label{eq:rc}
R(u, G) & = \frac{|V|!}{\prod_{w\in V} T^u_w}, 
\end{align}
%%%%%%
where $T^u_w$ is the size of the subtree of $G$ that is rooted at
$w$ and points away from $u$. For example, in Figure \ref{fig:rcexample}, 
let $u$ be node $1$. Then $|V| = 5$; the subtree sizes are 
$T^1_1 = 5$, $T^1_2 = 3$, $T^1_3= T^1_4 = T^1_5 = 1$ and
hence $R(1,G) = 8$. In \cite{ref:rc}, a linear time algorithm
is described to compute the rumor centrality of all nodes building on the
relation $R(u,G)/R(v,G) = T_u^v/T_v^u$ for neighboring nodes $u, v \in V$ ($(u,v) \in E$).

The rumor centrality of a given node $u \in V$ for a tree given by \eqref{eq:rc} 
is precisely the number of distinct spreading orders that could lead to the 
rumor infected graph $G$ starting from $u$. This is equivalent to computing 
the number of linear extensions of the partial order imposed by the graph 
$G$ due to causality constraints of rumor spreading. Under the SI model with homogeneous
exponential spreading times and a regular tree, it turns out that each of the spreading
orders is equally likely. Therefore, rumor centrality turns out to be the 
maximum likelihood (ML) estimator for the source in this specific setting (cf. \cite{ref:rc}). 
In general, the likelihood of each node $u \in V$ being the source given $G$ 
is proportional to the weighted summation of the number of distinct spreading 
orders starting from $u$, where weight of a spreading order could depend on 
the details of the graph structure and spreading time distribution of the SI
model. Now for a tree graph and SI model with homogeneous exponential spreading times, 
as mentioned above, such a quantity can be computed in linear time. But in general,
this could be complicated. For example, computing the number of linear extensions of
a given partial order is known to be \#P-complete (\cite{BrightwellWinkler}). While
there are algorithms for approximately sampling linear extensions given a partial
order (\cite{KK}), \cite{ref:rc} proposed the following simpler alternative for general graphs. 
%%%%%% 
\begin{dfn}\label{dfn:rc}[Rumor Centrality] Given node $u \in V$ in graph $G=(V,E)$, 
let $T \subset G$ denote a breadth-first search tree of $u$ with respect
to $G$. Then, the rumor centrality of $u$ with respect to $G$ is obtained by 
computing it as per \eqref{eq:rc} with respect to $T$. The estimated rumor
source is the one with maximal rumor centrality (ties broken uniformly at random). 
\end{dfn}
%%%

%%%%%%%%%%%%%%%%%%%%%%%%%%%%%%%%%%%%%%%%%%%%%%%%%%%%%%%%
%%%%%%%%%%%%%%%%%%%%%%%%%%%%%%%%%%%%%%%%%%%%%%%%%%%%%%%%%

\newcommand{\Pb}{\mathbf{P}}

%%%%%%%%%%%%%%%%%%%%%%%%%%%%%%%%%%%%5
\subsection{Prior Results} 

In \cite{ref:rc}, the authors established that rumor centrality is the maximum-likelihood estimator for the rumor 
source when the underlying graph $\G$ is a regular tree. They studied the effectiveness of this ML 
estimator for such regular trees. Specifically, suppose we observe the $n(t)$ node rumor infected graph $G$ after 
time $t$, which is a subgraph of $\G$. Let $C_{t}^k$ be the event that the source estimated as per rumor 
centrality is the $k$th infected node, and thus $C_{t}^1$ corresponds to the event of correct detection.  
The following are key results from \cite{ref:rc}: 
%%%%%%%%%
\begin{thm}[\cite{ref:rc}]\label{thm:rc1} Let $\G$ be a $d$-regular infinite tree with $d \geq 2$. Let  
\begin{align}\label{eq:rc1} 
\alpha_d^L = \liminf_{t\to\infty} \Pb\Big(C^1_{t}\Big) \leq \limsup_{t\to\infty} \Pb\Big(C^1_{t}\Big) = \alpha_d^U. 
\end{align}
Then, 
\begin{align}\label{eq:rc2}
\alpha_2^L = \alpha_2^U = 0, \quad
~\alpha_3^L = \alpha_3^U = \frac{1}{4}, \quad\mbox{and} \quad
0 < \alpha_d^L & \leq \alpha_d^U \leq \frac{1}{2}, ~\forall~d\geq 4. 
\end{align}
\end{thm}
%%%%%%%%%%%%%%%%%%%%%%%%%%%%%%%%%%%%%%%%%%

%%%%%%%%%%%%%%%%%%%%%%%%%%%%%%%%%%%%%%%%%%%%%%%%%%%%%%%%%%%%%%%%%%%%%%%%%%%%%%%%%%%%%%
%%%%%%%%%%%%%%%%%%%%%%%%%%%%%%%%%%%%%%%%%%%%%%%%%%%%%%%%%%%%%%%%%%%%%%%%%%%%%%%%%%%%%%
%%%%%%%%%%%%%%%%%%%%%%%%%%%%%%%%%%%%%%%%%%%%%%%%%%%%%%%%%%%%%%%%%%%%%%%%%%%%%%%%%%%%%%
%%%%%%%%%%%%%%%%%%%%%%%%%%%%%%%%%%%%%%%%%%%%%%%%%%%%%%%%%%%%%%%%%%%%%%%%%%%%%%%%%%%%%%
%%%%%%%%%%%%%%%%%%%%%%%%%%%%%%%%%%%%%%%%%%%%%%%%%%%%%%%%%%%%%%%%%%%%%%%%%%%%%%%%%%%%%%
%%%%%%%%%%%%%%%%%%%%%%%%%%%%%%%%%%%%%%%%%%%%%%%%%%%%%%%%%%%%%%%%%%%%%%%%%%%%%%%%%%%%%%
%%%%%%%%%%%%%%%%%%%%%%%%%%%%%%%%%%%%%%%%%%%%%%%%%%%%%%%%%%%%%%%%%%%%%%%%%%%%%%%%%%%%%%
%%%%%%%%%%%%%%%%%%%%%%%%%%%%%%%%%%%%%%%%%%%%%%%%%%%%%%%%%%%%%%%%%%%%%%%%%%%%%%%%%%%%%%

\section{Main Results} \label{sec:results}

We state the main results of this paper. In a nutshell, our results concern the characterization of 
the probability of $C^k_{t}$ for any $k \geq 1$ for large $t$ when $\G$ is a generic tree. 
As a consequence, it provides a characterization of the performance for sparse random
graphs.  

%%%%%%%%%%%%%%%%%%%%%%%%%%%%%%%%%%%%%%

\subsection{Regular Trees, SI Model with Exponential Spreading Times} %%%%%%%%%%%%%%%%%%%%%%%%%%%%%%%%%%%%%%%%%%%%%%%%%%%%%%%%%%%%

\begin{figure}\centering
	\centering
		\includegraphics[scale=0.5]{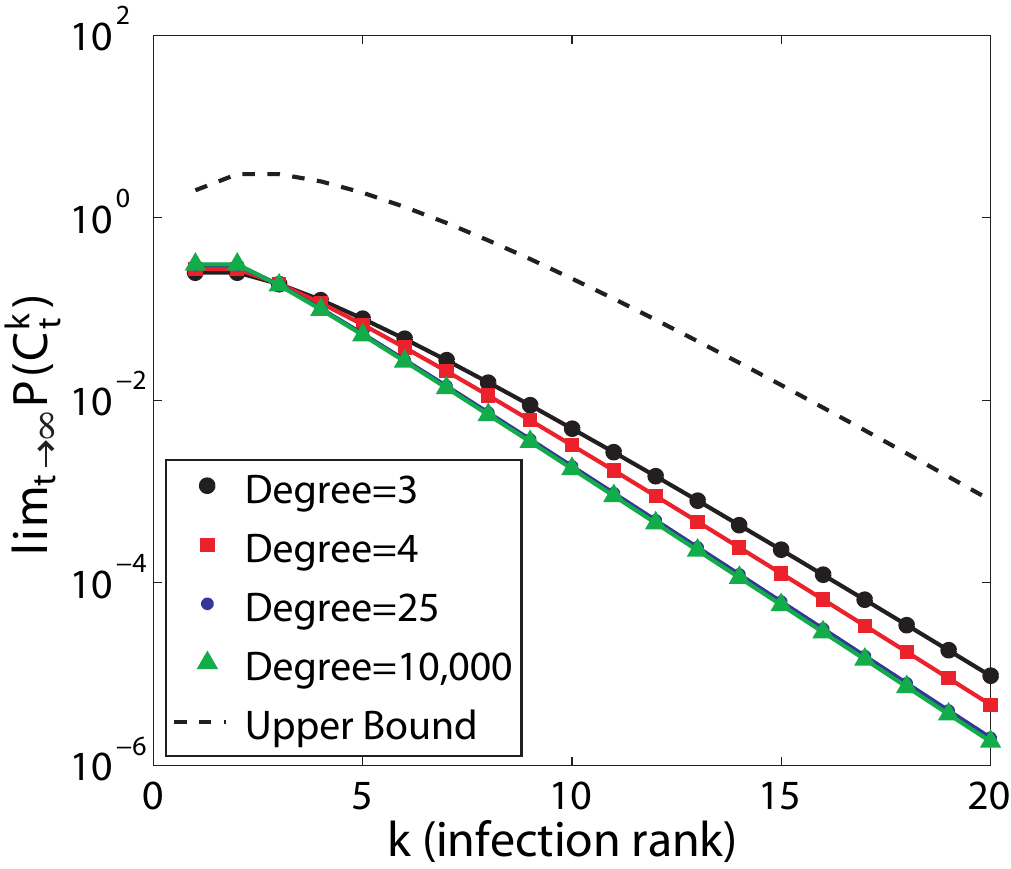}
	\caption{$\lim_{t\to\infty}\mathbf P\paranth{C^k_{t}}$ versus $k$ for regular trees of different degree.}
	\label{fig:regTreeError}
\end{figure}

\begin{figure}\centering
	\centering
		\includegraphics[scale=1]{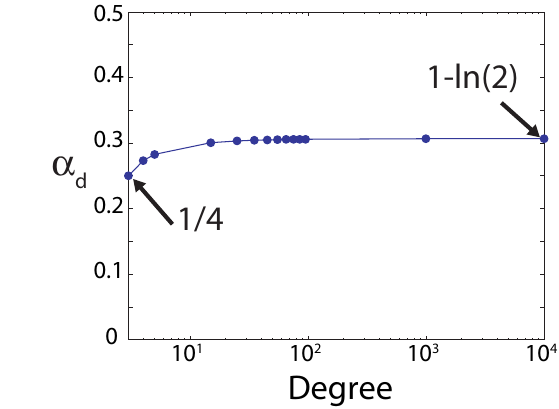}
	\caption{$\alpha_d$ versus degree $d$ for regular trees.  }
	\label{fig:regTreePd}
\end{figure}

We first look at rumor source detection on regular trees with degree $d\geq 3$, 
where rumor centrality is an exact ML estimator when the spreading times are exponentially distributed.  
Our results will utilize properties of Beta random variables.  We recall that the regularized incomplete 
Beta function $I_x(a,b)$ is the probability that a Beta random variable with parameters $a$ and $b$ is
less than $x \in [0,1]$, 
%%%%%%%%%%
\begin{align}\label{eq:beta}
I_x(a,b) & = \frac{\Gamma(a+b)}{\Gamma(a)\Gamma(b)} \int_0^x t^{a-1} (1-t)^{b-1} dt,
\end{align} 
%%%%%%%
where $\Gamma(\cdot)$ is the standard Gamma function.  For regular trees of degree $\geq 3$ we obtain the following result.
%%%%%%%%%%%%%%%%55
%%%%%%%%%%%%
\begin{thm}\label{thm:regTree}
Let $\G$ be $d$-regular infinite tree with $d \geq 3$. Assume a rumor spreads on $\G$ as per the SI model with exponential distribution with rate $\lambda$.  Then, for any $k \geq 1$, 
\begin{align}\label{eq:thmregTree}
\lim_{t\to\infty} \Pb\Big(C_{t}^k\Big) =&I_{1/2}\left(k-1+\frac{1}{d-2},1+\frac{1}{d-2}\right)\nonumber\\
        &+
	     (d-1)\left( I_{1/2}\left(\frac{1}{d-2},k+\frac{1}{d-2}\right) -1 \right). 
\end{align}
\end{thm}
%%%%%%%%%%%%%
%%%%%%%%%%%%%%%%%%%%%%%%%%%%%%%%%%%%%%%%%%%%%%%%%%%%%%%%%%%%%%%%%%%%%%%%%%%%%%%%%%%%%%%%%%%%%%%%%%%%%%%%%%%%%%%%%
For $k = 1$, Theorem \ref{thm:regTree} yields that 
$\alpha_d^L =\alpha_d^U = \alpha_d$ for all $d \geq 3$ where
\begin{align}\label{eq:alphad}
 \alpha_d & = d I_{1/2}\left(\frac{1}{d-2},\frac{d-1}{d-2}\right) - (d-1).
\end{align}
More interestingly, 
%%%%%%%%%% 
\begin{cor}\label{corollary:randRegInf}
\begin{align}\label{eq:ln2}
\lim_{d\to\infty} \alpha_d & = 1-\ln 2 ~\approx~0.307. 
\end{align}
\end{cor}
%%%%%%%%%%%%
%%%%%%%%%%%%%%%%%%%%%%%%%%%%%%%%%%%%%%%%%%%%%%%%%%%%%%%%%%%%%%%%%%%%%%%%%%%%%%%%%%%%%%%%
%%%%%%%%%%%%%%%%%%%%%%%%%%%%%%%%%%%%%%%%%%%%%%%%%%%%%%%%%%%%%%%%%%%%%%%%%%%%%%%%%%%%%%%%
For any $d \geq 3$, we can obtain a simple upper bound for 
Theorem \ref{thm:regTree} which provides the insight that 
the probability of error in the estimation decays exponentially 
with error distance (not number of hops in graph, but based on chronological
order of infection) from the true source. 
%%%%%%%%%%%%%%%%%%%%%%%%%%%%%%%%%%%%%%%%%%%%%%%%%%%%%%%%%%%%%%%%%%%%%%%%%%%%%%%%%%%%%%%
\begin{cor}\label{corollary:randReg3}
When $\G$ is a $d$-regular infinite tree, for any $k \geq 1$, 
\begin{equation*}
	\lim_{t\rightarrow \infty} \mathbf P\Big(C_{t}^k\Big)\leq  k\paranth{k+1}\paranth{\frac{1}{2}}^{k-1}~\asymp~\exp\Big(-\Theta(k)\Big).
\end{equation*}
\end{cor}
%%%%%%%%%%%%%%%%%%%%%%%%%%%%%%%%%%%%%%%%%
%%%%%%%%%%%%%%%%%%%%%%%%%%%%%
%%%%%%%%%%%%%%%%%%%%%%%%%%%%%%%%%%%%%%%%%%%%%%%%%%%%%%%%%%%%%%%%%%%%%%%%%%%%%%%%%%%%%%%%
To provide intuition, we plot the asymptotic error distribution $\lim_{t\to\infty}\mathbf P\paranth{C^k_{t}}$ for 
different degree regular trees in Figure \ref{fig:regTreeError}.  As can be seen, for degrees greater than $4$, all 
the error distributions fall on top of each other, and the probability of detecting the $k^{th}$ infection as 
the source decays exponentially in $k$.  We also plot the upper bound from Corollary \ref{corollary:randReg3}.  
As can be seen, this upper bound captures the rate of decay of the error probability.  Thus we see tight 
concentration of the error for this class of graphs. Figure \ref{fig:regTreePd} plots the asymptotic 
correct detection probability $\alpha_d$ versus degree $d$ for these regular trees.  It can be seen that 
the detection probability starts at $1/4$ for degree 3 and rapidly converges to $1-\ln(2)$ as the degree 
goes to infinity.

%%%%%%%%%%%%%%%%%%%%%%%%%%%%%%%%%%%%%%%%%%%%%%%%%%%%%%%%
%%%%%%%%%%%%%%%%%%%%%%%%%%%%%%%%%%%%%%%%%%%%%%%%%%%%%%%%%
%%%%%%%%%%%%%%%%%%%%%%%%%%%%%%%%%%%%%%%%%%%%%%%%%%%%%%%%%
%%%%%%%%%%%%%%%%%%%%%%%%%%%%%%%%%%%%%%%%%%%%%%%%%%%%%%%%%
%%%%%%%%%%%%%%%%%%%%%%%%%%%%%%%%%%%%%%%%%%%%%%%%%%%%%%%%%
%%%%%%%%%%%%%%%%%%%%%%%%%%%%%%%%%%%%%%%%%%%%%%%%%%%%%%%%
%%%%%%%%%%%%%%%%%%%%%%%%%%%%%%%%%%%%%%%%%%%%%%%%%%%%%%%%%

%%%%%%%%%%%%%%%%%%%%%%%%%%%%%%%%%%%%%%%%%%%%%%%%%%%%%%%%%%%%%%%%%%%%%%%%%%%%%%%%%%%%%%%%
%%%%%%%%%%%%%%%%%%%%%%%%%%%%%%%%%%%%%%%%%%%%%%%%%%%%%%%%%%%%%%%%%%%%%%%%%%%%%%%%%%%%%%%%
%%%%%%%%%%%%%%%%%%%%%%%%%%%%%%%%%%%%%%%%%%%%%%%%%%%%%%%%%%%%%%%%%%%%%%%%%%%%%%%%%%%%%%%%
%%%%%%%%%%%%%%%%%%%%%%%%%%%%%%%%%%%%%%%%%%%%%%%%%%%%%%%%%%%%%%%%%%%%%%%%%%%%%%%%%%%%%%%%

\subsection{Generic Random Trees, SI Model with Generic Spreading Times}\label{ssec:randTree} 

The above precise results were obtained using the memoryless property of the exponential distribution 
and the regularity of the trees.  Next, we wish to look at a more general setting both in terms of
tree structures and spreading time distributions. In this more general setting, while we cannot 
obtain precise values for the detection and error probabilities, we are able to make statements about 
the non-triviality of the detection probability of rumor centrality. When restricted to exponential
spreading times for generic trees, we can identify bounds on the error probability as well. Let us
start by defining what we mean by generic random trees through a generative model. 

\begin{dfn}[Generic Random Trees]
It is a rooted random tree, generated as follows: given a root node as a starting
vertex, add $\eta_0$ children to root where $\eta_0$ is an independent random variable 
with distribution $\D_0$. If $\eta_0 \neq 0$, then add a random number of children chosen 
as per distribution $\D$ over $\{0,1,\dots\}$ independently to each child of the root.  Recursively, 
to each newly added node, add independently a random number of nodes as per distribution $\D$. 
\end{dfn}

The generative model described above is precisely the standard Galton-Watson branching process
if $\D_0 = \D$. If we take $\D_0$ and $\D$ to be deterministic distributions with support on $d$ and
$d-1$ respectively, then it gives the $d$-regular tree. For a random $d$-regular graph on $n$ nodes,
as $n$ grows the neighborhood of a randomly chosen node in the graph converges (in distribution, locally) 
to such a $d$-regular tree. If we take $\D_0 = \D$ as a Poisson distribution
with mean $c > 0$, then it asymptotically equals (in distribution) to the local neighborhood of a 
randomly chosen node in a sparse Erdos-Renyi graph as the number of nodes grows.  Recall that a 
(sparse) Erdos-Renyi graph on $n$ nodes with parameter $c$ is generated by selecting each of 
the ${n\choose 2}$ edges to be present with probability $c/n$ independently. Effectively, random 
trees as described above capture the local structure for sparse random graphs reasonably well. 
For that reason, establishing the effectiveness of rumor centrality for source detection for such
trees provide insights into its effectiveness for sparse random graph models. 

We shall consider spreading time distributions to be generic. Let $F: [0,\infty) \to [0,1]$ be the 
cumulative distribution function of the spreading times. Clearly $F(0) = 0$, $F$ is non-decreasing and
$\lim_{t\to \infty} F(t) = 1$.
In addition, we shall require that the distribution is non-atomic at $0$, i.e. $F(0^+) = 0$.
We state the following result about the effectiveness of rumor centrality with such generic spreading 
time distribution. 

%%%%%%%%%%%%%%%%%%%%%%%%%%%%%%%%%%%%%%%%%%%%
\begin{thm}\label{thm:randtree}
Let $\eta_0$, distributed as $\D_0$, be such that $\Pr(\eta_0 \geq 3) > 0$ and 
let $\eta$, distributed as per $\D$, be such that $1<\Eb[\eta] < \infty$. Suppose the rumor starts from the root of the random tree generated as per
distributions $\D_0$ and $\D$ as described above and spreads as per the SI model with a 
spreading time distribution with an absolutely continuous density. Then, 
\[
\liminf_{t\to\infty} \Pb\Big(C^1_{t}\Big) > 0.  
\]
\end{thm}
%%%%%%%%%%%%%%%%%%%%%%%%%%%%%%%%%%%%%%%%%%5
The above result says that irrespective of the structure of the random trees, spreading time
distribution and  elapsed time, there is non-trivial probability of detecting the root as the 
source by rumor centrality. The interesting aspect of the result is that this non-trivial
detection probability is established by studying events when the tree grows without bound. 
{For finite size trees with $n$ nodes, the rumor source can be estimated by selecting
a random node, giving a probability of correct detection of $n^{-1}>0$.}  However, such
events are trivial and are not of much interest to us (neither mathematically, nor 
motivationally).

%%%%%%%%%%%%%%%%%%%%%%%%%%%%%%%%%%%%%%%%%%%%%%%%%%%%%%%%%%%%%%%%%%%%%%%%%%%

\subsubsection{Generic Random Trees, SI Model With Exponential Spreading Times} \label{ssec:randTree-exp} 

Extending the results of Theorem \ref{thm:randtree} for explicitly bounding the probability of the error event 
$\Pb\big(C^k_{t}\big)$ for generic spreading time distribution seems rather challenging. Here we
provide a result for generic random trees with exponential spreading times. 

%%%%%%%%%%%%%%%%%%%%%%%%%%%%%%%%%%%%%%%%%%%%%%%%%%%%%%%%%%%%%%%%%%%%%%%%%
%{\color{red} edit this.}
\begin{thm}\label{thm:randTreeError}
Consider the setup of Theorem \ref{thm:randtree} with spreading times being homogeneous exponential distributions 
with (unknown, but fixed) parameter $\lambda > 0$. In addition, let $\D_0 = \D$. Let $\eta$, 
distributed as per $\D$, be such that $\Eb[\eta] > 1$ and $\Eb[\exp(\theta \eta)] < \infty$ for
all $\theta \in (-\varepsilon, \varepsilon)$ for some $\varepsilon > 0$. 
Then, for appropriate constants $C', C'' > 0$,
\begin{align}
	 \limsup_{t\to\infty} \Pb\Big(C^k_{t}\Big) & \leq  C' \exp(-k C'').
	\end{align}
\end{thm}
%%%%%%%%%%%%%%%%%%%%%%%%%%%%%%%%%%%%%%%%%%%%%%%%%%%%%%%%%%%%%%%%%%%%%%%%%

The above result establishes an explicit upper bound on the {probability of the error event. The bound applies to essentially any generic random tree and demonstrates that
the probability of identifying later infected nodes as the rumor source decreases exponentially fast.}

%%%%%%%%%%%%%%%%%%%%%%%%%%%%%%%%%%%%%%%%%%%%%%%%%%%%%%%%%%%%%%%%%%%%%%%%%%%%%%
\subsubsection{Geometric Trees, SI Model With Generic Spreading Times}\label{ssec:geom}  
The trees considered thus far, $d$-regular trees with $d \geq 3$ or random trees with $\Eb[\eta] > 1$, 
grow exponentially in size with the diameter of the tree. This is in contrast with {path graphs} or $d$-regular
trees with $d = 2$ which grow only linearly in diameter. It can be easily seen that the probability of
correct detection, $\Pb(C^1_{t})$ will scale as $\Theta(1/\sqrt{t})$ for {path graphs} as long as the
spreading time distribution has non-trivial variance (see \cite{ref:rc} for proof of this statement for {the} SI model
with exponential spreading times). In contrast, the results of this paper stated thus 
far suggest that the expanding trees allow for non-trivial detection as $t \to \infty$. Thus, qualitatively 
path graphs and expanding trees are quite different -- one does not allow detection while the other does. To understand where the precise detectability threshold lies, we look at { polynomially
growing  {\em geometric} trees. }

\begin{dfn}[Geometric Tree] A geometric tree is a rooted, non-regular tree parameterized 
by constants $\alpha$, $b$, and $c$, with $\alpha\geq0$, $0 < b \leq c$, and root node $v^*$.  Let $d^*$ be
the degree of $v^*$, let the neighbors of $v^*$ be denoted
$v_1,v_2,...v_{d^*}$, and let the subtree rooted at $v_i$ and directed
away from $v^*$ be denoted
by  $T_i$ for $i=1,2,..., d^*$.  Denote  
 the number of nodes in $T_i$ at distance exactly $r$ 
from the subtree's root node $v_i$ as $n^i(r)$. Then we require that for all $1\leq i\leq d^*$ 
\begin{equation}
	b r^\alpha \leq n^i(r) \leq c r^{\alpha}. \label{eq:geom}
\end{equation}
\end{dfn}
The condition imposed by \eqref{eq:geom} states that each 
of the neighboring subtrees of the root should satisfy 
polynomial growth (with exponent $\alpha > 0$) and regularity 
properties. The parameter $\alpha > 0$ characterizes the 
growth of the subtrees and the ratio $c/b$ describes the 
regularity of the subtrees.  If $c/b\approx 1$ then the 
subtrees are somewhat regular, whereas if the ratio is much 
greater than 1, there is substantial heterogeneity in the subtrees.  
Note that the {path graph} is a geometric tree with $\alpha = 0$, 
$b=1$, and $c = 2$. 

We shall consider the scenario where the rumor starts from the root node
of a rooted geometric tree. We shall show that rumor centrality detects
the root as the source with an asymptotic probability of $1$ for a generic 
spreading time distribution with exponential tails. This is quite interesting
given the fact that rumor centrality is an ML estimator only for regular trees
with exponential spreading times. The precise result is stated next. 

%%%%%%%%%%%%%%%%%%%%%%%%%%%%%%%%%%%%%%%%%%%%%%%%%%%%%%%%%%%%%%%%%%%%%%%%%%%%%%%%%%%%%%%%%%%%%%
\begin{thm}\label{thm:geom}
Let $\G$ be a rooted geometric tree as described above with parameters $\alpha > 0$, $0 < b \leq c$ 
and root node $v^*$ with degree $d^*$ such that 
%%%%%%%%%
\begin{align*}
d_{v^*} & >  \frac{c}{b} + 1.  
\end{align*}
%%%%%%%%%%%%%%%
Suppose the rumor starts spreading on $\G$ starting from $v^*$ as per the SI model 
with a generic spreading time distribution whose cumulative density function $F: \Real \to [0,1]$
is such that (a) $F(0) = 0$, ~(b) $F(0^+) = 0$, and (c) if $X$ is a random variable distributed
as per $F$ then $\Eb[\exp(\theta X)] < \infty$ for $\theta \in (-\varepsilon, \varepsilon)$ for some $\varepsilon > 0$. Then
\begin{equation*}
	\lim_{t}\mathbf{P}(C_{t}^1)=1.
\end{equation*}
\end{thm}
A similar theorem was proven in \cite{ref:rc}, but only for the SI 
model with exponential spreading times.  We have now extended this result 
to arbitrarily distributed spreading times.
Theorem \ref{thm:geom} says that $\alpha=0$ and $\alpha>0$ serve as a threshold 
for non-trivial detection:  for $\alpha=0$, the graph is a 
path graph, so we would expect the detection probability to go to $0$ 
as $t\to \infty$ as discussed above, but for $\alpha>0$  the detection 
probability converges to $1$ as $t\to\infty$. 

%%%%%%%%%%%%%%%%%%%%%%%%%%%%%%%%%%%%%%%%%%%%%%%%%%%%%%%%%%5
\subsection{Detection Probability and Graph Growth: Discussion}\label{ssec:growth}
{Our results can be viewed as relating detection probability to graph growth parametrized by
$\alpha$.  For {path graphs}, where no detection is possible, $\alpha=0$.  For any finite, positive $\alpha$
we have geometric graphs where the detection probability converges to one.  For regular trees
or random graphs, the growth is exponential, which gives $\alpha=\infty$, 
and we have a detection probability that is strictly
between zero and one.}

{To understand these results at a high level, it is helpful to consider the properties of the rumor center given by Lemma \ref{lem:rc}.  Essentially, this lemma states that the graph is \emph{balanced} around the rumor center.  For the rumor source to be the rumor center (and therefore correctly identified as the true source), the rumor must spread in a balanced way.  For a {path graph} ($\alpha=0$), balance is a very delicate condition, requiring both subtrees of the source to be exactly equal in size.  The probability of this occurring goes to zero as the graph size goes to infinity.  }

{For any non-negative, finite alpha, this balance condition becomes easier to achieve if the source has degree greater than or equal to three.  In this case, because the number of vertices grows polynomially, the variation of the size of a rumor infected subtree after a time $t$ is much smaller than the expected value of its size, resulting in a concentration of the size.  This means that with high probability, no subtree will be larger than half of the network size, and balance is achieved.  The key here is that the boundary where the rumor can spread grows slower than the size of the rumor infected graph.  If the graph has $d^{\alpha+1}$ nodes, then the boundary contains $d^{\alpha}$ nodes.}

{For infinite alpha, which corresponds to graphs with exponential growth, the rumor boundary size is of the same order of magnitude in size as the rumor infected graph.  This results in a high variance in the subtree size.  We would expect this high variance to result in detection becoming impossible.  However, our analysis shows that the manner in which the rumor spreads on these graphs results in detection being possible with strictly positive probability.  Another way to view this result is that the vertices in each subtree act as witnesses which we can use to triangulate the source.  If there are three or more subtrees, and the subtree sizes do not vary considerably (as in graphs with polynomial growth), then the witnesses have low noise, and we can detect the source exactly as the observed rumor infected graph grows.  For exponentially growing graphs, the noise in the signals provided by the witnesses grows with the number of witnesses.  The increased number of witnesses balances the increased noise to give a detection probability that remains strictly positive as the graph size goes to infinity.}

%%%%%%%%%%%%%%%%%%%%%%%%%%%%%%%%%%%%%%%%%%%%%%%%%%%%%%%%%%%%%%%%%%%%%%%%%%%%%%%%%%%%%%%%%%%
%%%%%%%%%%%%%%%%%%%%%%%%%%%%%%%%%%%%%%%%%%%%%%%%%%%%%%%%%%%%%%%%%%%%%%%%%%%%%%%%%%%%%%%%%%%
\subsection{Locally Tree-Like Graphs: Discussion}\label{ssec:randgraph}

The results of the paper are primarily for tree structured graphs. On one hand, 
these are specialized graphs. On the other hand, they serve as local approximations for a variety of sparse
random graph models. As discussed earlier, for a random $d$-regular graph over $m$ nodes, a 
randomly chosen node's local neighborhood (say up to distance $o(\log m)$) is a tree with high
probability. Similarly, consider an Erdos-Renyi graph over $m$ nodes with each edge being present 
with probability $p=c/m$ independently for any $c > 0$ ($c > 1$ is an interesting regime due to the existence the of a giant component). A randomly chosen node's local neighborhood
(up to distance $o(\log m)$) is a tree and distributionally equivalent (in the large $m$ limit) to a random tree with Poisson degree distribution. 

Given such `locally tree-like' structural properties, if a rumor spreads on a random $d$-regular graph 
or sparse Erdos-Renyi graph for time $o(\log m)$ starting from a random node, then rumor 
centrality can detect the source with guarantees given by Theorems \ref{thm:regTree} and 
\ref{thm:randtree}. Thus,  although the results of this paper are for tree structured
graphs, they do have meaningful implications for tree-like sparse graphs.

\begin{figure}\centering
	\centering
		\includegraphics[scale=0.5]{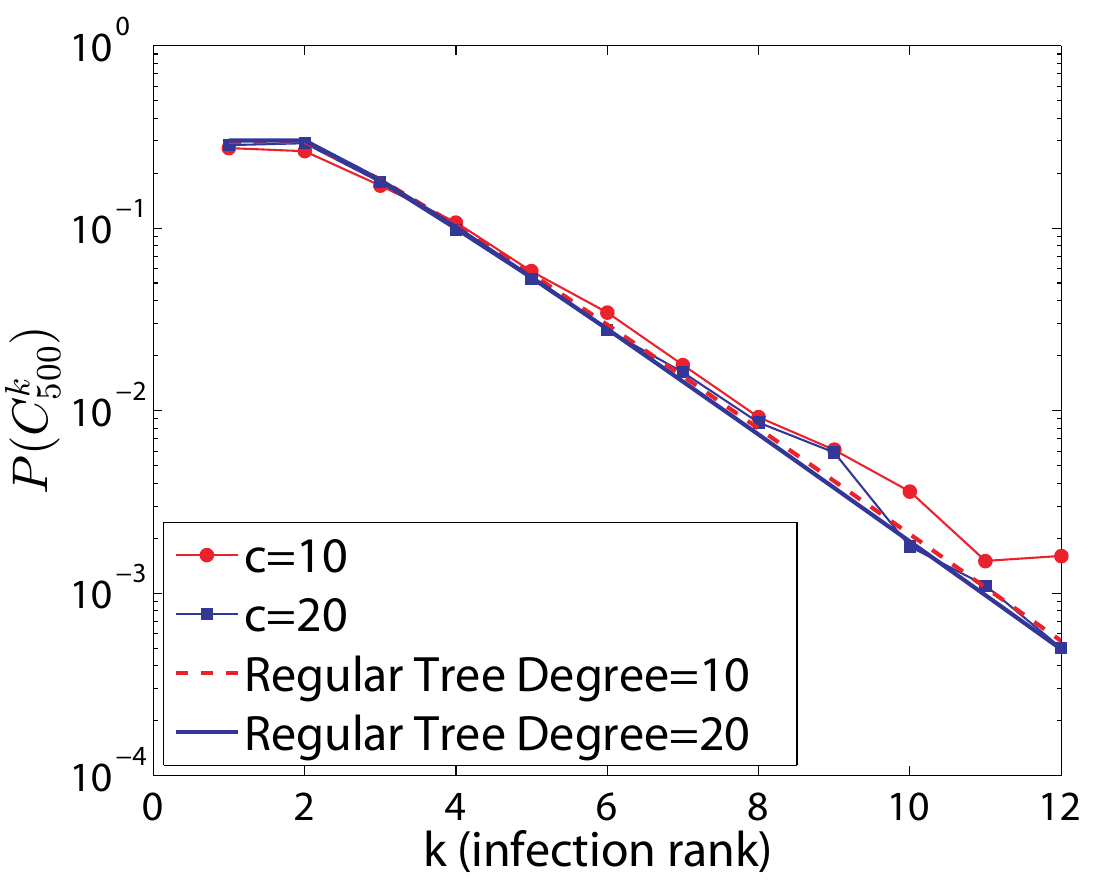}
	\caption{ {Empirical error probability versus $k$ 
	for Erdos-Renyi graphs with 500 vertices, mean degree $10$ and $20$, 
	and	exponentially distributed spreading times with mean one.  
	Also shown are $\lim_{t\to\infty} 
	\mathbf P\paranth{C^k_{t}}$ for degree $10$ and $20$ 
regular trees with exponentially distributed spreading times with mean one. } }
	\label{fig:ER}
\end{figure}

For the purpose of illustration, we conducted some simulations for Erdos-Renyi graphs that are reported in Figure \ref{fig:ER}. We generated
graphs with $m=50,000$ nodes and edge probabilities $p=c/m$ for $c=10$ and $c=20$.  
{The rumor graph contained $n=500$ nodes and the spreading times had an exponential distribution with mean one.  We used the general graph version of
rumor centrality as defined in Definition \ref{dfn:rc} as the rumor source estimator.}   We 
ran $10,000$ rumor spreading simulations to obtain the empirical error distributions plotted in Figure \ref{fig:ER}.   As can be seen, the 
error drops of exponentially in $k$, very similar to the regular tree error distribution. 
{ To make this more evident, we also plot the 
asymptotic error distributions for regular 
trees of degree $10$ and $20$ and it can be seen that the error decays at similar, exponential rates. } This indicates that even though 
there is substantial randomness in the graph, the asymptotic rumor source detection error distribution behaves as though it were 
a regular tree graph.  This result also suggests that the bounds in Theorem \ref{thm:randTreeError} are loose for this graph.

%%%%%%%%%%%%%%%%%%%%%%%%%%%%%%%%%%%%%%%%%%%%%%%%%%%%%%%%
%%%%%%%%%%%%%%%%%%%%%%%%%%%%%%%%%%%%%%%%%%%%%%%%%%%%%%%%%

%%%%%%%%%%%%%%%%%%%%%%%%%%%%%%%%%%%%%%%%%%%%%%%%%%%%%%%%%%%%%%%%%%%%%
%%%%%%%%%%%%%%%%%%%%%%%%%%%%%%%%%%%%%%%%%%%%%%%%%%%%%%%%%%%%%%%%%%%%%
%%%%%%%%%%%%%%%%%%%%%%%%%%%%%%%%%%%%%%%%%%%%%%%%%%%%%%%%%%%%%%%%%%%%%
%%%%%%%%%%%%%%%%%%%%%%%%%%%%%%%%%%%%%%%%%%%%%%%%%%%%%%%%%%%%%%%%%%%%%
%%%%%%%%%%%%%%%%%%%%%%%%%%%%%%%%%%%%%%%%%%%%%%%%%%%%%%%%%%%%%%%%%%%%%
\section{Proofs}\label{sec:proofs}
Here proofs of the results stated in Section \ref{sec:results} are presented.  We establish results for $d$-regular trees by connecting 
rumor spreading with Polya urn models and branching processes.  Later we extend this novel method to establish results for 
generic random trees under arbitrary spreading time distributions. After this, we prove Theorem \ref{thm:geom} using standard 
Chernoff's bound and the polynomial growth property of geometric trees.

%%%%%%%%%%%%%%%%%%%%%%%%%%%%%%%%%%%%%%%%%%%%%%%%%%%%%%%%%%%%%%%%%%%%%
\subsection{Proof of Theorem \ref{thm:regTree}: $d$-Regular Trees} 
%%%%%%%%%%%%%%%%%%%%%%%%%%
%%%%%%%%%%%%%%%%%%%%%%%%%%%%%%%%%%%%%%%%%%%%%%%%%%%%%%%%%%%%%
\subsubsection{Background: Polya's Urn.}\label{sec:polyaurn} We will recall Polya's urn process and it's asymptotic properties
that we shall crucially utilize in establishing Theorem \ref{thm:regTree}.  An interested reader can find a good exposition 
in \cite{ref:an}. 

{In the simplest form, Polya's Urn process operates in discrete time. Initially, at time $0$, an urn contains balls of two types, say 
$W_0$ white balls and $B_0$ black balls. Let $W_n and B_n$ denote the number of white and black balls, respectively, 
at the end of time $n \geq 1$. At each time $n\geq 1$, a ball is drawn at random from the urn 
($W_{n-1} + B_{n-1}$ balls in total). This ball is added back along with $\alpha \geq 1$ new balls of the same type leading 
to a new configuration of balls $(W_n, B_n)$. For instance, at time $n$, a white ball is drawn with probability $W_{n-1}/(W_{n-1} + B_{n-1})$  
and we have that $W_n = W_{n-1} + \alpha, ~B_n = B_{n-1}$.}

{Under the above described process, it is easy to check that the fraction of white (or black) balls is a bounded martingale. Therefore,
by the martingale convergence theorem, it has a limit almost surely. What is interesting is that the limiting distribution is nicely 
characterized as stated below. }
\begin{thm}{\cite[Theorem 1, pp. 220]{ref:an}}
For the Polya's Urn process described above
\begin{align}
	\frac{W_n}{W_n+B_n}& \rightarrow Y \qquad \mbox{almost surely},
\end{align}
where $Y$ is a Beta random variable with parameters $W_0/\alpha$ and $B_0/\alpha$. That is, 
for  $x \in [0,1]$,  
\begin{align*}
P(Y \leq x) & = I_x\Big(\frac{W_0}{\alpha}, \frac{B_0}{\alpha}\Big),
\end{align*}
where $I_x(a, b)$ is the incomplete Beta function defined as in \eqref{eq:beta}
\end{thm}

\subsubsection{Setup and Notation.}  
Let $\G=(\V,\E)$ be an infinite $d$-regular tree and let the rumor start 
spreading from a  node, say $v_1$. Without loss of generality, we view the tree as a randomly generated tree, as described in Section \ref{sec:results}, with $v_1$ being the root with 
$d$ children and all the subsequent nodes with $d-1$ children (hence each node has degree $d$). 
We shall be interested in $d \geq 3$.  Now suppose the rumor is spread on this tree starting 
from $v_1$ as per the SI model with exponential distribution with rate $\lambda>0$. 

Initially, node $v_1$ is the only rumor infected node and its $d$ neighbors are potential nodes 
that can receive the rumor.  We will denote the set of nodes that are not yet rumor infected but are
neighbors of rumor infected nodes as the {\em rumor boundary}. Initially the rumor
boundary consists of the $d$ neighbors of $v_1$. Under the SI model, each edge 
has an  independent exponential clock of mean $1/\lambda$.  The minimum of $d$ 
independent exponentials of mean $1/\lambda$ is an exponential random variable of mean $1/\paranth{d\lambda}$,
and hence when one of the $d$ nodes (chosen uniformly at random) in the rumor boundary gets 
infected, the infection time has an exponential distribution with mean $1/\paranth{d\lambda}$. Upon this infection, 
this node gets removed from the boundary and adds its $d-1$ children to the rumor boundary. 
That is, each infection adds $d-2$ new nodes to the rumor boundary. 
In summary,  let $Z(t)$ denote the number of nodes in the rumor boundary at time $t$, 
then $Z(0)=d$ and $Z(t)$ evolves as follows: each of the $Z(t)$ nodes has an 
exponential clock of mean $1/\lambda$; when it ticks, it dies and {$d-1$ 
new nodes are born} which in turn start their own independent exponential clocks of mean $1/\lambda$
and so on. 
%{The thus described $Z(t)$ is precisely a Markov branching process (MBP).
%If we think of the contributions of each of the $d$ sub-trees of
%the root $v_1$ to the rumor boundary separately we effectively have $d$ 
%branching processes each starting at time $0$ with initial value equal to $1$}. 
 Let 
$u_1,\dots, u_d$ be the children of $v_1$; let $Z_i(t)$ denote the number of nodes in the
rumor boundary that belong to the {subtree $\mathcal T_i(t)$} that is rooted at $u_i$ with 
$Z_i(0) = 1$  for $1\leq i\leq d$; $Z(t) = \sum_{i=1}^d Z_i(t)$.  {Let $T_i(t)=|\mathcal T_i(t)|$} 
denote the total number of nodes infected in the subtree 
rooted at $u_i$ at time $t$; initially $T_i(0) = 0$ for $1\leq i\leq d$. Since each infected
node add $d-2$ nodes to the rumor boundary, it can be easily checked that 
$Z_i(t) = (d-2) T_i(t) + 1$ and hence $Z(t) = (d-2) T(t) + d$ with $T(t)$ being the
total number of infected nodes at time $t$ (excluding $v_1$). 

%

%%%%%%%%%%%%%%%%%%%%%%%%%%%%%%%%%%%%%%%%%%%%%%%%%%
%%%%%%%%%%%%%%%%%%%%%%%%%%%%%%%%%%%%%%%%%%%%%%%%%%
%%%%%%%%%%%%%%%%%%%%%%%%%%%%%%%%%%%%%%%%%%%%%%%%%%
\subsubsection{Probability of Correct Detection. } 
Suppose we observe the rumor infected nodes at some time
$t$ which we do not know. That is, we observe the rumor infected graph $G(t)$ which contains
the root $v_1$ and its $d$ infected subtrees $\mathcal T_i(t)$ for $1\leq i\leq d$.  We recall the following
result of \cite{ref:rc} that characterizes the rumor center {(for a proof see Section 
\ref{sec:proof_lem_rc}). }

\begin{lem}[\cite{ref:rc}]\label{lem:rc} Given a {tree} graph $G=(V,E)$, there can be at most
two rumor centers. Specifically, a node $v \in V$ is  a rumor center if and only if 
\begin{align}\label{eq:rc-char}
T_i^v & \leq \frac{1}{2}\Big(1+\sum_{j\in \N(v)} T_j^v\Big), ~~~~\forall ~i \in \N(v),
\end{align}
where $\N(v) = \{u \in V: (u,v) \in E\}$ are neighbors of $v$ in $G$ and $T^v_j$ 
denotes the size of the sub-tree of $G$ that is rooted at node $j \in \N(v)$ and includes all nodes that are away from node $v$ (i.e. the subtree does not include $v$). The 
rumor center is unique if the inequality in \eqref{eq:rc-char} is strict for all $i \in \N(v)$.
\end{lem}
%%%%

This immediately suggests the characterization of the event that node $v_1$, the
true source, is identified by rumor centrality at time $t$: $v_1$ is
a rumor center only {if $2T_i(t) \leq 1+\sum_{j=1}^d T_j(t)$ for} all $1\leq i\leq d$, 
and if the inequality is strict then it is the unique rumor center. 
 {Let $E_i = \{2T_i(t) < 1+\sum_{j=1}^d T_j(t)\}$
and $F_i = \{2T_i(t) \leq 1+\sum_{j=1}^d T_j(t)\}$}. Then, 
\begin{align}
\Pb\Big(C^1_{t}\Big) & \geq \Pb\Big(\cap_{i=1}^d E_i\Big) %\nonumber \\
%                                       & 
~= 1- \Pb\Big(\cup_{i=1}^d E_i^c\Big) %\nonumber \\
%                                       & 
~\stackrel{(a)}{\geq} 1- \sum_{i=1}^d \Pb\Big(E_i^c\Big) %\nonumber \\
                                 %      & 
~\stackrel{(b)}{=} 1 - d \Pb\Big(E_1^c\Big). \label{eq:dreg1}
\end{align} 
Above, (a) follows from the union bound of events and (b) from symmetry. Similarly, 
we have 
\begin{align}
\Pb\Big(C^1_{t}\Big) & \leq \Pb\Big(\cap_{i=1}^d F_i\Big) %\nonumber \\
                                       % & 
~= 1- \Pb\Big(\cup_{i=1}^d F_i^c\Big)% \nonumber \\
                                       ~\stackrel{(a)}{=} 1- \sum_{i=1}^d \Pb\Big(F_i^c\Big)% \nonumber \\
                                       ~ \stackrel{(b)}{=} 1 - d \Pb\Big(F_1^c\Big). \label{eq:dreg2}
\end{align} 
Above, (a) follows because events $F_1^c,\dots, F_d^c$ are disjoint and (b) from
symmetry. Therefore, the probability of correct detection boils down to evaluating
$\Pb(E_1^c)$ and $\Pb(F_1^c)$ which, as we shall see, will coincide with each 
other as $t\to\infty$. Therefore, the bounds of \eqref{eq:dreg1}
and \eqref{eq:dreg2} will provide the exact evaluation of the correct detection 
probability as $t\to\infty$.
%%%%%%%%%%%%%%%%%%%%%%%%%%%%%%%%%%%%%%%%%%%%%%%%%%
%%%%%%%%%%%%%%%%%%%%%%%%%%%%%%%%%%%%%%%%%%%%%%%%%%

\subsubsection{$\Pb(E_1^c)$, $\Pb(F_1^c)$ and Polya's Urn.} 
Effectively, the interest is in {the ratio $T_1(t)/(1+\sum_{i=1}^d T_i(t))$, especially} as 
$t \to \infty$ (implicitly we are assuming that this ratio is well defined for a given $t$ or 
else by definition there is only one node infected which will be $v_1$, the true source). 
It can be easily verified that as $t\to \infty$, $T_i(t) \to \infty$ for all $i$ almost surely and 
hence $Z_i(t) = (d-2) T_i(t) + 1$ goes to $\infty$ as well. Therefore, it is sufficient to study 
the ratio $Z_1(t)/(\sum_{j=1}^d Z_j(t))$ as $t \to \infty$ since we shall find that this ratio 
converges to a random variable with density on $[0,1]$. In summary, if we establish that 
the ratio $Z_1(t)/(\sum_{j=1}^d Z_j(t))$ converges in distribution on $[0,1]$ with a well 
defined density, then it immediately follows that $\Pb(E_1^c) \stackrel{t\to\infty}{\longrightarrow} \Pb(F_1^c)$  
and we can use  $Z_1(t)/(\sum_{j=1}^d Z_j(t))$ in {place of   $T_1(t)/(1+\sum_{j=1}^d T_j(t))$.}

{With these facts in mind, let us study the ratio $Z_1(t)/(\sum_{j=1}^d Z_j(t))$. 
For this, it is instructive to view the simultaneous evolution of $(Z_1(t), Z_{\neq 1}(t))$ 
(where $Z_{\neq 1}(t) \stackrel{\triangle}{=} \sum_{j=2}^d Z_j(t)$) as that  induced by the 
standard, discrete time, Polya's urn. Initially, $\tau_0 =0$ and there is one ball of type $1$ (white)
representing $Z_1(\tau_0) = 1$ and $d-1$ balls of type $2$ (black) representing $Z_{\neq 1}(\tau_0) = d-1$
in a given {\em urn}. The $j^{th}$ event happens at time $\tau_j$ (also known as a split time) 
when one of the $Z_1(\tau_{j-1}) + Z_{\neq 1}(\tau_{j-1})$ ($= d+(j-1)(d-2)$)  
{balls chosen uniformly  at random is returned to the urn along with $d-2$ new balls of its type. }
If we set $\tau_j - \tau_{j-1}$ equal to an exponential random variable with mean 
$1/(\lambda(d + (j-1) (d-2)))$, then it is easy to check that the fraction of balls of type $1$ 
is identical in law to that of $Z_1(t)/(\sum_{i=1}^d Z_i(t))$ (here we are using
the {\em memoryless} property of exponential random variables crucially). 
Therefore, for our purposes, it is sufficient to study the limit law of fraction of 
balls of type $1$ (or white) under this Polya's urn model. }

%It is well understood that the fraction of balls of type $1$ at time $\tau_j$, 
%which is equal to $Z_1(\tau_j)/(Z_1(\tau_j) + Z_{\neq 1}(\tau_j))$, is a 
%martingale with value in $[0,1]$. By the standard martingale convergence 
%theorem, it converges to a well defined random variable almost surely. 
%The limiting distribution of this random variable is given by 
%the following theorem.
%\begin{thm}{\cite[Theorem 1, pp. 220]{ref:an}}
%Consider a Polya's urn that begins with $W_0\geq 0$ balls of type 1
%and $B_0\geq 0$ balls of type 2 and upon each draw of a ball,
%returns the ball and adds $\alpha\geq 1 $ balls of the same type.
%Let $W_n$ and $B_n$ be the number of type 1 and type 2 balls, respectively,
%after $n$ draws.  Then 
%\begin{align}
%	\frac{W_n}{W_n+B_n}& \rightarrow Y \qquad \mbox{almost surely},
%\end{align}
%where $Y$ is a Beta random variable
%with parameters $W_0/\alpha$ and $B_0/\alpha$.
%\end{thm}
From the discussion in Section \ref{sec:polyaurn}, it follows that the ratio 
$Z_1(t)/(\sum_{i=1}^d Z_i(t))$ converges to a Beta random variable with 
parameters  $1/(d-2)$ and $(d-1)/(d-2)$. Since the Beta distribution has a density
on $[0,1]$, from the above discussion it follows that as  $t \to \infty$, 
$|\Pb(E_1^c) - \Pb(F_1^c)| \to 0$ and hence from \eqref{eq:dreg1}, \eqref{eq:dreg2}
%%%%%%%%%%%
\begin{align}
\lim_{t\to\infty} \Pb\Big(C^1_{t}\Big) & = 1 - d \Big(1-I_{1/2}\Big(\frac{1}{d-2}, 1+\frac{1}{d-2}\Big)\Big),
\end{align}
%%%%%%%%%%%%%%%
where $I_{1/2}(a,b)$ is the probability that a Beta random variable 
with parameters $a$ and $b$ takes value in $[0,1/2]$. Note that this establishes
the result of Theorem \ref{thm:regTree} for $k = 1$ in \eqref{eq:thmregTree}. 
%%%%%%%%%%%%%%%%%%%%%%%%%%%%%%%%%%%%%%%%%%%%%%%%%%%
%%%%%%%%%%%%%%%%%%%%%%%%%%%%%%%%%%%%%%%%%%%%%%%%%%%%

\subsubsection{Probability of $C^k_{t}$.} \label{sssec:Cnk}
Thus far we have established Theorem 
\ref{thm:regTree} for $k = 1$ (the probability of the rumor center being the true source).
The probability of the event $C_{t}^k$ (the $k$th infected node being the rumor center)
is evaluated in an almost identical manner with a minor difference. For this reason,
we present an abridged version of the proof. 

%%%%%%%%%%%%%%%%%%%%%%%%%%%%%%%%%%%%%%%%%%%%%%%%%%%%%%%%
\begin{figure}\centering
	\includegraphics[scale=1]{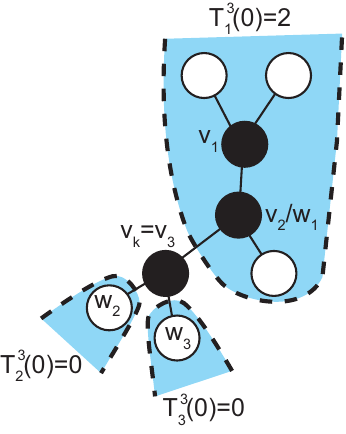}
	\caption{Illustration of the labeling of the neighbors of $v_k$
	and their subtrees for $k=3$ in a rumor graph at $t_k$ 
	(the time of infection of $v_k$).  
	The rumor infected nodes are colored black, 
	and the uninfected nodes are white.}
	\label{fig:cnkneighbors}
\end{figure}
%%%%%%%%%%%%%%%%%%%%%%%%%%%%%%%%%%%%%%%%%%%%%%%%%%%%%%%%%

Let $T_k = \inf\{t : n(t) = k\}$ represent the time when $k^{th}$ node is infected. It can be easily checked that for $d \geq 3$ regular
tree with exponential spreading time distribution, $T_k < \infty$ with probability $1$. Consider $t \geq T_k$. For $k \geq 2$, let $v_k$  
be the $k^{th}$ infected node when the rumor starts from $v_1$. We will evaluate the probability of identifying $v_k$ as the 
rumor center. Let $G$ represent the rumor infected tree observed at time $t$ with $n(t) \geq k$ nodes. Let $w_1,\dots, w_d$ be 
the $d$ neighbors of $v_k$, as is illustrated in Figure \ref{fig:cnkneighbors}.  We shall denote the 
neighbor of $v_k$ that is along the path joining $v_k$ and $v_1$ as $w_1$. Note that $w_1$ must have been infected before
$v_k$ when the rumor starts spreading from $v_1$. Let $w_2,\dots, w_{d}$ be the $d-1$ `children' of $v_k$, away from $v_1$.  

For convenience, we shall use notation $t' = t - T_k$ with $t' \geq 0$.  Let $\mathcal T_i^k(t')$ be the subtree 
of $G$ rooted at $w_i$ at time $t$ away from $v_k$.  Therefore, $\mathcal T_1^k(t')$ is rooted at $w_1$ and includes 
$v_1,\dots, v_{k-1}$. For $ ~2\leq i\leq d$,  $\mathcal T_i^k(t')$ are rooted at $w_i$ and contain nodes in $G$ that are 
away from $v_k$.  None of the $\mathcal T_i^k(t')$ for $1\leq i\leq d$
include $v_k$.  When $v_k$ is infected at time $T_k$, we have that $T^k_1(0) = k-1$, and $T^k_i(0)=0$ for $2\leq i\leq d$.  
This notation is illustrated in Figure \ref{fig:cnkneighbors}.   

By definition $\mathcal T_1^k(t')$ 
is never empty, but $\mathcal T_i^k(t')$  can be empty if 
$w_i$ is not infected, for $2\leq i\leq d$. As before, let 
$T_i^k(t')=|\mathcal T_i^k(t')|$. As per Lemma  \ref{lem:rc}, $v_k$ 
is identified as a rumor center if and only if all of 
its $d$ subtrees are balanced, i.e. 
%%%%%%%%%%%
\begin{align}
2T_i^k(t') & \leq 1+\sum_{j=1}^d T_j^k(t'), ~~\forall~1\leq i\leq d.
\end{align}
%%%%%%%%%%%
Therefore, for $t \geq T_k$ with $t' = t - T_k$, 
 \begin{align}
\Pb\Big(C^k_{t}\Big) & \geq \Pb\Big(\cap_{i=1}^d E_i\Big)~= 1- \Pb\Big(\cup_{i=1}^d E_i^c\Big) {\geq} 1- \sum_{i=1}^d \Pb\Big(E_i^c\Big), \quad \mbox{and} \label{eq:dreg3}\\ 
\Pb\Big(C^k_{t}\Big) & \leq \Pb\Big(\cap_{i=1}^d F_i\Big) ~= 1- \Pb\Big(\cup_{i=1}^d F_i^c\Big)~{=} 1- \sum_{i=1}^d \Pb\Big(F_i^c\Big),\label{eq:dreg4}
\end{align} 
{where $E_i = \{2T^k_i(t') < 1+\sum_{j=1}^d T^k_j(t')\}$ 
and $F_i = \{2T^k_i(t') \leq 1+\sum_{j=1}^d T^k_j(t')\}$}. 

As before, we shall evaluate these probabilities by studying the evolution of the appropriate rumor boundaries. However, unlike for $k = 1$, 
when $k \geq 2$ the rumor boundaries have asymmetric initial conditions. Specifically,
$\mathcal T_1^k(0) = k-1$, $Z^k_1(\cdot)  = (d-2)(k-1) +1$ and for $2\leq i\leq d$, 
$\mathcal T_i^k(0) = 0$ and $Z_i^k(0) = 1$. Beyond this difference, the rules governing the evolution of the
rumor boundaries are the same as those described in the proof for $k = 1$. To evaluate $E_1^c$ (and $F_1^c$), we consider 
a Polya's urn in which we start with $ (d-2)(k-1) +1$ balls
of type $1$ (corresponding to $Z_1^k(0)$) and $d-1$ balls of type $2$ (corresponding to $\sum_{j=2}^d Z_j^k(0)$). 
With these initial conditions, the limit law of fraction of balls of type $1$ turns out to be (see \cite{ref:an} for details) 
a Beta distribution with parameters $a = ((d-2)(k-1)+1)/(d-2)~=~(k-1) + 1/(d-2)$ and $b = (d-1)/(d-2) = 1 + 1/(d-2)$. 
Finally, since the fraction of balls of type $1$, i.e. the ratio $Z_1^k(t')/(\sum_{j=1}^k Z_j^k(t'))$, equals { $T_1^k(t')/(1+\sum_{j=1}^d T_j^k(t'))$} as $t'\to\infty$, we obtain 
%%%%%%%%%%%%
\begin{align}\label{eq:dreg5}
\lim_{t'\to\infty} \Pb\Big(E_1^c\Big) ~=~\lim_{t'\to\infty} \Pb\Big(F^c_1\Big)~=~ 1 - I_{1/2}\Big(k-1 + \frac{1}{d-2}, ~1+ \frac{1}{d-2}\Big). 
\end{align}
%%%%%%%%%%%%%%%%%
For $2\leq i\leq d$, in the corresponding Polya's urn model, we start
with $1$ ball of type $1$ and $k (d-2) + 1$ balls of type $2$.  Therefore, using an identical sequence of arguments,
we obtain that for $2\leq i\leq d$, 
%%%%%%%%%%%%%%%%%%%%%
\begin{align}\label{eq:dreg6}
\lim_{t'\to\infty} \Pb\Big(E_i^c\Big) ~=~\lim_{t'\to\infty} \Pb\Big(F^c_i\Big)~=~ 1 - I_{1/2}\Big(\frac{1}{d-2}, ~k+ \frac{1}{d-2}\Big). 
\end{align}
%%%%%%%%%%
From \eqref{eq:dreg3}-\eqref{eq:dreg6}, it follows that 
%%%%%%%%%%%%%%%
\begin{align}
\lim_{t\to\infty} \Pb\Big(C^k_{t}\Big)  = &I_{1/2}\Big(k-1 + \frac{1}{d-2}, ~1+ \frac{1}{d-2}\Big)
 + (d-1)\Big(I_{1/2}\Big(\frac{1}{d-2}, ~k+ \frac{1}{d-2}\Big)-1\Big).
\end{align}
%%%%%%%%%%%
This establishes \eqref{eq:thmregTree} for all $k$ and completes the proof of Theorem \ref{thm:regTree}. 
%%%%%%%%%%%%%%%%%%%%%%%%%%%%%%%%%%%%%%%%%%%%%%%%%%%%%%%%%%
%%%%%%%%%%%%%%%%%%%%%%%%%%%%%%%%%%%%%%%%%%%%%%%%%%%%%%%%%%
%%%%%%%%%%%%%%%%%%%%%%%%%%%%%%%%%%%%%%%%%%%%%%%%%%%%%%%%%%

\subsection{Proof of Lemma \ref{lem:rc}}\label{sec:proof_lem_rc}
{We provide here a proof of Lemma \ref{lem:rc} for the convenience
of the reader.  Much of this proof is taken from \cite{ref:rc}.
We begin by establishing the following property about rumor centrality.}
\begin{prop}\label{prop:rc1}
Consider an undirected tree graph $G=(V,E)$ with $|V|=N$ 
and any two neighboring nodes $u,v\in V$ such
that $(u,v)\in E$.  The rumor centralities of 
these two nodes satisfy the following
relationship:

\begin{align}
	         \frac{R(u,G)}{R(v,G)} =\frac{T^v_u}{N-T^v_u}.
\end{align}
\end{prop}
%%%%%%%%%%%%%%%%%%%%%%%%%%%%%%%%%%%%%%%%%%%%%%%%%%%%%
We now show that if $v$ is a rumor center then it must satisfy the condition given by equation \eqref{eq:rc-char} in Lemma \ref{lem:rc}.  For any node $i$ neighboring the rumor center $v$, Proposition \ref{prop:rc1} gives
\begin{align*}
	\frac{R(i,G)}{R(v,G)}&=\frac{T^{v}_{i}}{N-T^{v}_{i}}\leq 1.
\end{align*}
Rearranging terms, we obtain
\begin{align*}
	T^v_i &\leq \frac{N}{2}%\\
	      %&
	      \leq \frac{1}{2}\paranth{1+\sum_{j\in\mathcal N(v)} T^v_j}.
\end{align*}
%%%%%%%%%%%%%%%%%%%%%%%%%%%%%%%%%%%%%%%%%%%%%%%%%%%
We now establish the other direction of Lemma \ref{lem:rc}.  
Assume equation \eqref{eq:rc-char} of the Lemma is satisfied for a node $v$.  We now
show that $v$ must be a rumor center.

 Let $i\in V$ be a node
$d$ hops from $v$ and let
$\curly{v_0=v,v_1,v_2,...,v_d=i}$ be the sequence of nodes in the path between $v$ and $i$.
  Using Proposition
\ref{prop:rc1} we obtain

\begin{align}
	    \frac{R(i,G)}{R(v,G)} & = \prod_{i=1}^d\frac{R(v_{i},G)}{R(v_{i-1},G)}%\nonumber\\
			%& 
			= \prod_{i=1}^d\frac{T^{v_{i-1}}_{v_i}}{N-T^{v_{i-1}}_{v_i}}\nonumber.
\end{align}
The subtrees on the path between $v$ and $i$ have the special property 
that $T^{v_{i-1}}_{v_i} = T^{v}_{v_i}$ for $i=1,2,...,d$ because the nodes in the subtree
rooted at $v_i$  are the same if the subtree is directed away from $v_{i-1}$ or $v$. 
We also have the property that $N/2\geq T^{v}_{v_{i-1}}> T^{v}_{v_i}$  for $i=2,...,d$ because
the subtrees must decrease in size by at least one node as we traverse the path from $v$ to $i$
and node $v$ satisfies equation \eqref{eq:rc-char} in the Lemma.  With these facts we obtain 
\begin{align}
	    \frac{R(i,G)}{R(v,G)} & = \prod_{i=1}^d\frac{T^{v}_{v_i}}{N-T^{v}_{v_i}}%\nonumber\\
			  %&
			  \leq 1.
\end{align}

If the inequality is strict in equation \eqref{eq:rc-char}, then we have that
for any $i\neq v$,
$T^v_i<N/2$.  Using Proposition \ref{prop:rc1} it can be shown that 
this implies that for every $i\neq v$, there exists
a node $j\neq i$ such that $T^i_j>N/2$.  This violates equation \eqref{eq:rc-char},
which means $i$ cannot be a rumor center.  Therefore, $v$ is the unique rumor center.
%%%%%%%%%%%%%%%%%%%%%%%%%%%%%%%%%%%%%%%%%%%%%%%%%%%%%%%%%%%%%%%%%

\subsection{Proof of Proposition \ref{prop:rc1}}

The rumor centrality of a node $v$ in a tree $G = (V,E)$ with $|V|=N$ is given by

\begin{align*}
   R(v,G) = \frac{N!}{\prod_{w\in V}T^v_w}
\end{align*}
with the tree variables $T^v_w$ denoting the size of the subtree of $G$ that is rooted
at $w$ and points away from $v$.  For any two nodes $u,v$ in a tree such that $(u,v)\in E$
there is a special relationship between their subtrees.  For any $w\in V, w\neq u,v$,
it can be shown that $T^v_w = T^u_w$.  Also, it can be shown that $\mathcal T^v_u$ contains
all nodes which are not in $\mathcal T^u_v$.  This gives the simple relation that $T^u_v = N-T^v_u$.
With these results on the subtree variables we obtain
\begin{align*}
   \frac{R(u,G)}{R(v,G)} &= \frac{\prod_{w\in V}T^v_w}{\prod_{w\in V}T^w_v}%\\
	%& 
	= \frac{T^v_u}{N-T^v_u}.
\end{align*}

%%%%%%%%%%%%%%%%%%%%%%%%%%%%%%%%%%%%%%%%%%%%%%%%%%%%%%%%%%%%%%%%%%%%%%%%%%%%%%%%
\newcommand{\deltaReg}{\frac{1}{d-2}}
\subsection{Proof of Corollary \ref{corollary:randRegInf}}
Simple analysis yields Corollary \ref{corollary:randRegInf}.  We start by defining the asymptotic 
probability for a $d$-regular tree as $\lim_{t\to\infty} \Pb\paranth{C_{t}^1}=\alpha_d$.  
This quantity then becomes

\begin{align*}
\alpha_d =&dI_{1/2}\left(\frac{1}{d-2},1+\frac{1}{d-2}\right) -d+1\\
         =&1-\frac{d\Gamma(1+\frac{2}{d-2})}{\Gamma(\deltaReg)\Gamma(1+\deltaReg)} \int_{\frac{1}{2}}^{1} t^{\deltaReg-1} (1-t)^{\deltaReg}dt\\
\end{align*}
We then take the limit as $d$ approaches infinity.
\begin{align*}
\lim_{d\rightarrow \infty} \alpha_d =&\lim_{d\rightarrow \infty} 1- \frac{d\Gamma(1+\frac{2}{d-2})}{\Gamma(\deltaReg)\Gamma(1+\deltaReg)} \int_{\frac{1}{2}}^{1} t^{\deltaReg-1} (1-t)^{\deltaReg}dt\\
          &= 1- \lim_{d\rightarrow \infty}  \frac{d\Gamma(1+\frac{2}{d-2})}{\paranth{d-2-\gamma +O\paranth{d^{-1}}}\Gamma(1+\deltaReg)} \int_{\frac{1}{2}}^{1} t^{\deltaReg-1}\paranth{1-t}^{\deltaReg}dt\\
          &= 1-\int_{\frac{1}{2}}^{1} t^{-1}dt\\
          &=1-\ln\paranth{2}.
\end{align*}
Above, $\gamma$ is the Euler-Mascheroni constant and we have used the following approximation of $\Gamma(x)$ for small $x$:	$\Gamma\paranth{x}= x^{-1}-\gamma +O\paranth{x}$.

%%%%%%%%%%%%%%%%%%%%%%%%%%%%%%%%%%%%%%%%%%%%%%%%%%%%%%%%%%%%%%%%%%%%%%%%%%%%%%%%

%%%%%%%%%%%%%%%%%%%%%%%%%%%%%%%%%%%%%%%%%%%%%%%%%%%%%%%%%%%%%%%%%%%%%%%%%%%%%%%%%%%%%%%%
\subsection{Proof of Corollary \ref{corollary:randReg3}}

\newcommand{\beps}{\varepsilon}

Corollary \ref{corollary:randReg3} follows from \eqref{eq:thmregTree} and monotonicity of the $\Gamma$ function 
over $[1,\infty)$.  For $k \geq 2$, 
\begin{align*}
	\lim_{t\rightarrow \infty} \mathbf P\Big(C_{t}^k\Big)  & =  I_{1/2}\left(k-1+\frac{1}{d-2},1+\frac{1}{d-2}\right)\\
    & \qquad \qquad + (d-1)\left(I_{1/2}\left(\frac{1}{d-2},k+\frac{1}{d-2}\right) - 1 \right) \\
& \leq I_{1/2}\left(k-1+\frac{1}{d-2},1+\frac{1}{d-2}\right)\\
	& = \frac{\Gamma(k+\frac{2}{d-2})}{\Gamma(k-1+\deltaReg)\Gamma(1+\deltaReg)} \int_0^{\frac{1}{2}} t^{k+\deltaReg-2} (1-t)^{\deltaReg} dt\\
          & \stackrel{(a)}{\leq} \frac{\Gamma(k+\frac{2}{d-2})}{\Gamma(k-1+\deltaReg)\Gamma(1+\deltaReg)}  \int_0^{\frac{1}{2}} t^{k-2}dt\\
	&\stackrel{(b)}{\leq} \frac{4e^2 \Gamma(k+2)}{\Gamma(k-1)}\int_0^{\frac{1}{2}} t^{k-2}dt\\
	& \stackrel{(c)}{\leq} 4e^2 k (k+1)(k+2) \paranth{\frac{1}{2}}^{k-1} \\
	& \asymp \exp\big( - \Theta(k)\big).
\end{align*}
%
%	&\stackrel{(b)}{\leq} \frac{4e^2 \Gamma(k+2)}{\Gamma(k-1)\Gamma(1)}\int_0^{\frac{1}{2}} t^{k-2}dt\\
%	&\leq \frac{4 e^2 k(k+1)}{k-1} \paranth{\frac{1}{2}}^{k-1}\\
%	& = \exp\Big(-k \ln 2 + \beps(k)\Big) ~=~\exp\big(-\Theta(k)\big),
In above, (a) follows from the fact that $t <1$ and hence $t^{k-2+1/(d-2)} \leq t^{k-2}$. 
For (b), we use the following well-known properties of the $\Gamma$ function: (i) over $[2,\infty)$, the $\Gamma$ function is non-decreasing and hence for $k\geq 2$ and $d\geq 3$, $\Gamma(k+\frac{2}{(d-2)}) \leq \Gamma(k+2)$; (ii) over $(0,\infty)$, the $\Gamma$
function achieves its minimal value in $[1,2]$ which is at least $\frac{1}{2e}$ and therefore, along with (i), we have that 
$\Gamma(k-1 + \frac{1}{d-2}) \geq \frac{\Gamma(k-1)}{2e}$ and $\Gamma(1+\frac{1}{d-2}) \geq \frac{1}{2e}$. For (c), we
use the fact that $\Gamma(x+1) = x \Gamma(x)$ for any $x \in (0,\infty)$. 

%we need to use property of the $\Gamma$ function carefully. It is well known that $\Gamma$ achieves minimum
%value in $[1,2]$, $\Gamma(1) = \Gamma(2) = 1$, and $\Gamma(x)$ increases for $x \geq 2$. Now the minimum value of 
%$\Gamma$ function is at least $1/(2e)$. Since $d \geq 3$, $1+1/(d-2) \in [1,2]$. Therefore, 
%$\Gamma(1+1/(d-2)) \leq \Gamma(1)/(2e)$. Similarly, $\Gamma(k-1 + 1/(d-2)) \geq \Gamma(k-1)/(2e)$ for all $k \geq 2$
%and $d \geq 3$. Therefore, (b) follows. 

%%%%%%%%%%%%%%%%%%%%%%%%%%%%%%%%%%%%%%%%%%%%%%%%%%%%%%%%%
%%%%%%%%%%%%%%%%%%%%%%%%%%%%%%%%%%%%%%%%%%%%%%%%%%%%%%%%%%%%%%%%%%%%%%%%%%%%%%%%%%%%%%%%%%%%%%%%%%%%%%%%%%%%%%%%%%%%%%%%%%%%%%%%%%%%%%%%%%%%%%%%%%%%%%%%%%%%%%
%%%%%%%%%%%%%%%%%%%%%%%%%%%%%%%%%%%%%%%%%%%%%%%%%%%%%%%%%%%%
\subsection{Proof of Theorem \ref{thm:randtree}: Correct Detection for Random Trees}\label{ssec:randtreeproof}

The goal is to establish that there is a strictly positive probability of detecting the
source correctly as the rumor center when the rumor starts at the root of a generic
random tree with generic spreading time distribution as defined earlier. The probability 
is with respect to the joint distribution induced by the tree construction and 
the SI rumor spreading  model with independent spreading times. We extend the
technique employed in the proof of Theorem \ref{thm:regTree}. However, it requires using
a generalized Polya's urn or age-dependent  branching process as well as 
delicate technical arguments. 

\subsubsection{Background: Age-Dependent  Branching Process.}\label{sec:gpu}
We recall a generalization of the classical Polya's urn known as an age-dependent
branching process. Such a process starts at time $t = 0$ with a given finite number of nodes, say $B(0) \geq 1$. 
Each node remains alive for an independent, identically distributed lifetime with cumulative distribution function given by $F: [0, \infty) \to [0,1]$. The lifetime distribution function $F$ will
be assumed to be non-atomic at $0$, i.e. $F(0^+) = 0$. Each node dies after remaining alive for its lifetime. Upon the death of a node, it gives birth to random number of nodes, say $\eta$. The random variables $\eta$ corresponding
to each node are independent and identically distributed over the non-negative integers. The newly born nodes live for their lifetime
and the upon death give birth to new nodes, and so on.  

As can be seen, the classical Galton-Watson process is a special case of this general model and the size of the entire urn in 
the Polya's urn process described earlier naturally fits this model. An interested reader is referred to \cite{ref:an} for a detailed 
exposition. Next we recall certain remarkable asymptotic properties of this process that will be crucially utilized. We start
with a useful definition. 
%%%%%%%
\begin{dfn}{\cite[pp. 146]{ref:an}}\label{dfn:malthusian}
Let $m \equiv {\mathbb E}[\eta]$.  The Malthusian parameter $\alpha = \alpha(m ,F)$  of an age-dependent branching process is
the unique solution, if it exists, of the equation 
\begin{align}
	m\int_{0}^{\infty}e^{-\alpha y}dF(y)=1.
\end{align}
\end{dfn}
%%
%%%%%%%%%%%%%%%%%%%%%%%%%
A sufficient condition for the existence of the Malthusian parameter is $m =  {\mathbb E}[\eta] > 1$. As an example, consider process
where spreading time distribution is exponential with parameter $\lambda$, i.e. $F(t) = 1-e^{-\lambda t}$, and let $m = \Eb[\eta] > 1$.  
The Malthusian parameter $\alpha(m,F)$ is given by the solution of 
%%%%%%%%%
\begin{align*}
	 m\int_{0}^{\infty}e^{-\alpha y}\lambda e^{-\lambda y}dy & =1,
\end{align*}
which is
\begin{align*}
	  \alpha(m,F) = \lambda (m-1).
\end{align*}
The Malthusian parameter captures the average growth rate of the branching process. We now recall the following result. 
%%%%%%%%%%%%%%%%%%%%%%%%%%%%%%%%%%%%%%%%%%%%%%%%%%%%%%%%%%%%%%%%%%%%
\begin{thm}{\cite[Theorem 2, pp. 172]{ref:an}}\label{thm:fact}
Consider an age-dependent branching process as described above with the additional properties that $m = \Eb[\eta] > 1$ and 
$\Eb[\eta \log \eta]  < \infty$. Let $\alpha \equiv \alpha(m, F)$ be the Malthusian parameter of the process and define 
\begin{align*}
	c & = \frac{m-1}{\alpha m^2\int_{0}^{\infty}ye^{-\alpha y}dF(y)}.
\end{align*}
Let $B(t)$ denote the number of nodes alive in the process at time $t \geq 0$. Then  
\begin{align*} 
\frac{1}{ce^{\alpha t}} B(t) \stackrel{t\to\infty}{\rightarrow} W \qquad \mbox{in~distribution}, 
\end{align*}
where $W$ is such that 
\begin{align}
\mathbf E[W]&=1\\
\Pb\big(W = 0\big) & = q, \label{eq:bp1} \\
\Pb\big( W \in (x_1, x_2)\big) & = \int_{x_1}^{x_2} w(y) dy, ~\quad \mbox{for}\quad 0 < x_1 < x_2 < \infty, \label{eq:bp2}
\end{align}
%%%%%%%%5
where $q \in (0,1)$ is the smallest root of the equation 
$\sum_{k=0}^\infty s^k \Pb(\eta=k) = s$ and $w(\cdot)$ 
is absolutely continuous with respect to the Lebesgue 
measure so that $\int_0^\infty w(y) dy = 1-q$. 
\end{thm}
%%%%%%%%%%%%%%%%%%%%%%%%%%%%%%%%%%%%%%%%%%%%%%%%%
{The above result states that with probability $q$ ($0 < q < 1$) the branching process becomes extinct, and with probability $1-q$
the size of the process scales as $\exp(\alpha t)$ for large $t$. We will need finer control on the asymptotic growth of the branching 
process. }Precisely, we shall use the following implication of the above stated result. 
%%%%%%%%%%%%%%%%%%%%%%%%%%%%%%%%%
\begin{cor}\label{cor:key}
Under the setting of Theorem \ref{thm:fact}, for any $f > 1$, there exists an $x > 0$ so that 
\begin{align}
\Pb\big(W \in (x, fx) \big) & >  0. 
\end{align}
\end{cor}
%
%\begin{proof}
\textbf{Proof}  Define 
\begin{align*}
a_k & = f^k, \quad \text{for}~~k \in {\mathbb{Z}}.
\end{align*}
%%%%%%%%%%%%%%%%%%%%%%%%
By definition, $\{W > 0 \} = \cup_{k \in \mathbb{Z}} ~\curly{W\in(a_k, a_{k+1}]}$. Due to the absolute continuity of $w(\cdot)$ in \eqref{eq:bp2}, it follows that $\Pb(W = a_k) = 0$ for all $k \in \mathbb{Z}$. Therefore, it follows 
that 
\begin{align}
0 & <  \Pb\big(W > 0\big) \nonumber \\
        & = \Pb\big(\cup_{k \in \mathbb{Z}} \big\{W \in (a_k, a_{k+1})\big\}\big) \nonumber \\
        & \leq \sum_{k \in \mathbb{Z}} \Pb\big(W \in (a_k, a_{k+1}\big). 
\end{align} 
From above, it follows that there exists a $k$ such that $\Pb\big(W \in (a_k, a_{k+1})\big) > 0$. This completes the proof.

%%%%%%%%%%%%%%%%%%%%%%%%%%%%%%%%%%%%%%%%%%%%%%%%%%%%%%%%%%%%%%%%%%%%%%%%%%%%%%%%%%%%%%

\subsubsection{Notation.} 
{We quickly recall some notation. To start with, as before let $v_1$ be the root node 
of the tree. It has $\eta_0$ children distributed as per $\D_0$. Define the event $A = \{\eta_0 \geq 3\}$.
By assumption of  Theorem \ref{thm:randtree}, $\Pb(A) > 0$. We shall show that }
\begin{align}
\liminf_{t\to\infty} \Pb(C^1_t | A) & > 0, 
\end{align}
{as it will imply the desired result $\liminf_{t\to\infty} \Pb(C^1_t) > 0$, using $\Pb(A) > 0$ since} $\Pb(C^1_t )  \geq \Pb(C^1_t | A ) \Pb(A)$. 
{Therefore, we shall consider conditioning on event $A$ and let $d = \eta_0 \geq 3$ for remainder of
the proof. Note that all the spreading times as well as all other randomness are independent of $\eta_0$.
The only effect of conditioning on $A$ is that we know that root has $d \geq 3$ children. 
Let $u_1,\dots, u_d$ be the $d$ children of root $v_1$. The random tree 
$\G$ is constructed by adding a random number of children to $u_1,\dots, u_d$ recursively 
as per distribution $\D$ as explained in Section \ref{ssec:randTree}. } 
%
%Note that the spreading
%times are independent of the degree of the root note and therefore, conditioning on the
%root having $d \geq 3$ children does not affect the spreading time process. Specifically, 
%the rumor starts at $v_1$ at time $0$ and spreads as per the SI model on $\G$ with spreading 
%times whose cumulative density function $F : \Real \to [0,1]$ is such that $F(0) = 0$, 
%$F(0^+) = 0$ and $\lim_{t\to\infty}F(t) = 1$.

{Further, as explained in Section \ref{ssec:randTree}, the rumor spreads on $\G$ starting from $v_1$ at
time $0$ as per the spreading times with cumulative distribution
function $F$ that is non-atomic at $0$. Let $G$ be the sub-tree of $\G$ that is infected at time $t$ with 
$n(t)$ infected nodes in $G$ at time $t$. Let $\mathcal T_i(t)$ denote the 
subtree of $G$ rooted at node $u_i$ (pointing away from root $v_1$) at time $t$, for $1\leq i\leq d$ and
let $T_i(t)=|\mathcal T_i(t)|$.  By definition $T_i(0) = 0$ for $1\leq i\leq d$. Let $Z_i(t)$ denote 
the size of  the rumor boundary of $\mathcal T_i(t)$; initially $Z_i(0) = 1$, ~$1\leq i\leq d$. }

%
%Since we are interested in evaluating the probability of detection with 
%respect to the joint distribution over the tree generation and the SI spreading
%model,  we model the evolution of $Z_i(\cdot)$ and $T_i(\cdot)$ as follows. 

{Now let us consider the evolution of $Z_i(\cdot)$: recall that each node in the rumor boundary 
has a rumor infected parent (neighbor). 
This node will become infected after the amount of time given by the spreading time associated
with the edge connecting the node with its infected parent. After the node becomes
infected, it is no longer part of the rumor boundary, but all of its uninfected neighbors (children) 
become part of the rumor boundary. And as per the random generative process of the 
tree construction, the number of children added, $\eta$, has distribution ${\cal D}$. 
Therefore, the rumor boundary process $Z_i(\cdot)$ for each $1\leq i\leq d$ is exactly an age-dependent
branching process. Further, each $Z_i(\cdot)$ evolves independently and since initially each starts at the same time
with exactly one node, they are identically distributed. Therefore, we can 
utilize the results stated in Section \ref{sec:gpu} to characterize the properties of 
$Z_i(\cdot)$ for $1\leq i\leq d$. In the case of regular trees, $Z_i(\cdot)$ and $T_i(\cdot)$ were linearly 
related which allowed us to obtain results about $T_i(\cdot)$ and the desired conclusion. 
While in this general setting, $Z_i(\cdot)$ and $T_i(\cdot)$ are not linearly related,  we show 
that they are asymptotically linearly related due to an appropriate Law of Large Numbers effect. This
will help us obtain the desired conclusion. We present the details next. }

\subsubsection{Correct Detection.} {As before, we wish to show that} 
\begin{align}
\Pb\Big(C^1_{t} | A \Big) & \geq \Pb\Big(\cap_{i=1}^d \Big\{2T_i(t) <1+ \sum_{j=1}^d T_j(t)\Big\}  \Big),  \label{eq:rt1} 
\end{align}
{where we have removed the conditioning on $A$, as the only effect of $A$ was having $d$ distinct trees, which is already captured. 
We shall establish \eqref{eq:rt1} in two steps:}
\begin{itemize}

\item[] {\em Step 1.} {Using the characterizations of $Z_i(\cdot)$ in terms of  age dependent branching processes as discussed 
above, we shall show that there is a non-trivial event  $\E_1 \subset \cap_{i=1}^d \Big\{2Z_i(t) < \sum_{j=1}^d Z_j(t)\Big\}$ with 
$\liminf_{t\to\infty} \Pb(\E_1) > 0$.}

\item[]{\em Step 2.} {Identify an event $\E_2 \subset \E_1$ with $\liminf_{t\to\infty} \Pb(\E_2) > 0$ and $\E_2 \subset  \cap_{i=1}^d \Big\{2T_i(t) < 1+\sum_{j=1}^d T_j(t)\Big\}$
for all $t$ large enough. }
\end{itemize}
{This will yield the desired results. }

%and then 
%establish a non-trivial lower bound on $\E$ using Theorem \ref{thm:fact}; ~(ii) establish that
%for $t$ large \textbf{enough, $\E \subset \cap_{i=1}^d \Big\{2T_i(t) < 1+\sum_{j=1}^d T_j(t)\Big\}$}. 
%This will immediately imply that $\Pb(\E)$ is a non-trivial lower bound on $\Pb\Big(C^1_{t}\Big)$.
%%%%%%%%%%%%%%%%%%%%%%%%%%%%%%%%%%%%%%%%%%%%%%%%%%%%%%%%%

\subsubsection{Step 1.} %A Non-Trivial Event} 
For any 
$x > 0$ and $\bepsilon > 0$ define the event $\E(x,\bepsilon, t)$ as 
\begin{align}
\E(x,\bepsilon, t) & = \cap_{i=1}^d \Big\{ Z_i(t) c^{-1}e^{-\alpha t} \in (x, (1-3\bepsilon) (d-1) x) \Big\}. 
\end{align}
Since $ d \geq 3$, $(1-3\bepsilon) (d-1) > 1$ for small enough $\bepsilon > 0$ and hence the above
event is well defined. It can be easily checked that 
$\E(x,\bepsilon, t)  \subset \cap_{i=1}^d \Big\{2Z_i(t) < \sum_{j=1}^d Z_j(t)\Big\}$, since under this
event, 
\[
\max_{i} Z_i(t) \leq  (1-3\bepsilon) (d-1) x < (d-1) x \leq (d-1) \min_{i} Z_i(t).
\]
By Theorem \ref{thm:fact}, it follows that $Z_i(t) c^{-1}e^{-\alpha t}$ converges to $W_i$, 
which are independent across $i$ and identically distributed as per \eqref{eq:bp1}-\eqref{eq:bp2}. 
Therefore, using Corollary \ref{cor:key}, it follows that  there exists an $x^* > 0$ such that 
\begin{align}\label{eq:rt3}
\liminf_{t\to\infty}\Pb\Big(\E(x^*,\bepsilon, t) \Big) & > 0.
\end{align}
Define $\E_1 = \E(x^*,\bepsilon, t)$. 

%%%%%%%%%%%%%%%%%%%%%%%%%%%%%%%%%%
%%%%%%%%%%%%%%%%%%%%%%%%%%%%%%%%%%%%%%%%%%%%%%%%%%%%%%%%%%%%%%

\subsubsection{Step 2.} 
{We want to find $\E_2 \subset \E_1$ so that for $t$ large enough, $\E_2 \subset \cap_{i=1}^d \Big\{2T_i(t) < 1+\sum_{j=1}^d T_j(t)\Big\}$
and $\liminf_{t\to\infty}  \Pb(\E_2) > 0$. For regular trees this was achieved by using the linear (deterministic) relationship
between the $Z_i(\cdot)$ and $T_i(\cdot)$. Here, we do not have such a relationship. Instead, we shall establish an asymptotic 
relationship.  To that end, recall that for any $t \geq 0$, }
\begin{align}
Z_i(t) & = 1 + \sum_{\ell \in \mathcal T_i(t)} (\eta_\ell -1).
\end{align}
{The above holds because as per the branching process, when a node in the `boundary' dies ($-1$ is added to $Z_i(\cdot)$) and 
it is added to $\Tin_i(\cdot)$,  $\eta_\ell$ new nodes are added to boundary. }

{Consider $\Tin_i(\cdot)$. It grows by adding nodes with a random number of children as per distribution $\D$ independently. Let $\eta_1,\eta_2,\dots$ be
these random number of children added to it in that order (we assume this sequence to be infinite irrespective of whether or not $\Tin_i(\cdot)$ stops growing). 
Since these are i.i.d. random variables with finite mean (actually, $\Eb[\eta \log \eta] < \infty$), by the standard Strong Law of Large Numbers, for 
any small enough $\bepsilon, \delta > 0$, with probability at least $1-\delta$, for all $1\leq i\leq d$, we have that for all $p \geq 1$}
%%%%%%%%%%%%%%%%%%%%%%
\begin{align}\label{eq:renew}
\frac{(1-\bepsilon) p}{m} - C(\bepsilon, \delta) & \leq N_i(p) ~\leq~ \frac{(1+\bepsilon)p}{m} + C(\bepsilon, \delta) 
\end{align}
%%%%%%%%%%%%%%%%%%%%%%%%%%%%%%%%
{where $N_i(p) = \inf\{\ell: \sum_{j=1}^\ell (\eta_j-1) \geq p\}$, $m = \Eb[\eta]$ and $C(\bepsilon, \delta)$ is a non-negative constant depending upon $\bepsilon, \delta$ but independent of $p$.  Let us call the event represented by \eqref{eq:renew} as $\E'(\bepsilon, \delta)$. Here,
we have the freedom of choosing as small a $\delta$ and $\bepsilon$ as we like. We will choose $\delta$ so that it is much smaller than the probability of
event $\E_1$ for $t$ large enough. Given such a choice, it will follow that for all $t$ large enough, 
the event $\E_2 = \E_1 \cap \E'(\bepsilon, \delta)$
has strictly positive probability. Under event $\E_2$, we have (with the definition $\hat{Z}_i(t) = Z_i(t)c^{-1}e^{-\alpha t}$) }
%%%%
\begin{align}\label{eq:end1}
\hat{Z}_i(t) & \in (x^*, x^*(1-3\bepsilon)(d-1)), ~~\text{for~all}~~1\leq i\leq d, \\
\Tin_i(t)c^{-1}e^{-\alpha t} & \in \paranth{\frac{\hat{Z}_i(t)(1-\bepsilon)}{m} - a_t, \frac{\hat{Z}_i(t)(1+ \bepsilon)}{m} + a_t}
 ~~\text{for~all}~1\leq i\leq d, \nonumber 
\end{align}
%%%
{where the constants $a_t \to 0$ as $t\to\infty$. Therefore, it can be easily checked that for $t$ large enough and $\bepsilon$ small enough, 
$\E_2 \subset \cap_{i=1}^d \Big\{2T_i(t) < 1+\sum_{j=1}^d T_j(t)\Big\}$, just the way we argued that $\E_1  \subset \cap_{i=1}^d \Big\{2Z_i(t) < \sum_{j=1}^d Z_j(t)\Big\}$. As discussed above, with an appropriate choice of $\delta$ and $\bepsilon$, we can guarantee that $\liminf_{t\to\infty} \Pb(\E_2) > 0$.  This concludes the search for the desired event $\E_2$ and we have established the desired claim of $\liminf_{t\to\infty} \Pb(C^1_{t}) > 0$. This completes the proof of Theorem \ref{thm:randtree}.}

%%
%\begin{align}\label{eq:end1}
%Z_i(t)c^{-1}e^{-\alpha t} & \in (x^*, x^*(1-3\bepsilon)(d-1)), ~~\text{for~all}~~1\leq i\leq d, \\
%\Tin_i(t)c^{-1}e^{-\alpha t} & \in \paranth{\frac{(Z_i(t)-1)(1-\bepsilon)}{m}, \frac{(Z_i(t)-1)(1+\bepsilon)}{m}}, ~~\text{for~all}~1\leq i\leq d. \nonumber 
%\end{align}
%%

%%%%

%%%%%%%%%%%%%%%%%%%%%%%%%%%%%%%%%%%%%%%%%%%%%%%%%%%%%%%%%%%%%%%%%%%%%%%%
%%%%%%%%%%%%%%%%%%%%%%%%%%%%%%%%%%%%%%%%%%%%%%%%%%%%%%%%%%%%%%%%%%%%%%%%
\subsection{Proof of Theorem \ref{thm:randTreeError}}

\subsubsection{Background: Properties of Age-Dependent Branching Processes.} { We shall utilize the following
property known in the literature about bounds on the moment generating function of the size of an age-dependent 
branching process. We shall assume the notation from the earlier section. }

\begin{thm}{\cite[Theorem 3.1]{nakayama2004finite}} \label{thm:mgf}
Consider an age dependent branching process with the properties that $m = \Eb[\eta] > 1$,  
$\Eb[\exp(\theta \eta)]  < \infty$ for all $\theta \in (0,\theta_1)$ for some $\theta_1 > 0$, 
and the spreading time distribution is non-atomic.  Let $B(t)$ represent the number of living nodes in the branching process at time $t$ and let $V(t)$ represent the number of nodes born before time $t$. Then, there exists a $\theta^*>0$ such that for all $\theta \in (-\theta^*, \theta^*)$
\begin{align}
\mathbf E\bracket{e^{\theta B(t)}} \leq \mathbf E\bracket{e^{\theta V(t)}}<\infty.
\end{align}
\end{thm}

\subsubsection{Background: Two inequalities.}

{We state two useful concentration-style inequalities that we shall derive here for completeness.}  
\begin{prop}\label{prop:std} For $i\geq 1$ let
$X_i$ be independent and identically distributed random variables  such 
that $\Eb[\exp(\theta X_1)] < \infty$ for all $\theta \in (-\delta, \delta)$ for some $\delta > 0$. 
Then, for any $\varepsilon>0$, there exists constants $C_1, C_2(\varepsilon, \delta) >0$ such that
\begin{align}
	\Pb\paranth{\sum_{i=1}^{n} X_i\leq \mu n(1-\varepsilon)} & \leq C_1 \exp\big(-C_2(\varepsilon, \delta) \mu n\big), 
\end{align}
where $\mu = \Eb[X_1]$. 
\end{prop}

\begin{prop}\label{prop:error}
Consider independent and identical random variables $X_1,\dots, X_{r+s}$ for integers 
$r, s$ such that $1 \leq s < r$. Let $\mu=\Eb[X_1]$ and $\Eb[\exp(\theta X_1)] < \infty$ for all $\theta \in (-\delta, \delta)$ 
for some $\delta > 0$. Then there exists a constant $c$ such that 
for any $\gamma > 0$, there exists a constant 
$\theta^* = \min(\frac{\gamma + (r-s)\mu}{2 (r+s) c}, \delta_1/2)$ for some $0 < \delta_1 < \delta$, such that 
\begin{align}
\Pb(\sum_{i=1}^r X_i - \sum_{j=1}^s X_{r+j} \leq -\gamma) & \leq  \exp\big(- \frac{1}{2}\theta^* (\gamma + (r-s) \mu)\big).
\end{align}
\end{prop}

{Next, we prove these two propositions. }

\medskip
\noindent
{\bf Proof of Proposition \ref{prop:std}.} {Let $X$ be a random variable with identical distribution as that of $X_i, i\geq 1$. 
By assumption in the Proposition statement, it follows that for $\theta \in (-\delta, \delta)$}
\begin{align}
M_X(\theta) & \equiv \log \Eb[\exp(\theta X)] \nonumber \\
 & = \log\Big(1 + \sum_{j=1}^\infty \theta^j \Eb[X^j]/j!\Big) \nonumber \\
 & \leq \log\Big(1 + \theta \mu + c \theta^2\Big),  \nonumber
\end{align}
%%%
{for some $c > 0$ for all $\theta \in (-\delta_1, \delta_1)$ for some $0 < \delta_1 < \delta$. Using the inequality $\log (1 + x) \leq x$ for all $x > -1$, we obtain }
\begin{align}\label{eqd:1}
M_X(\theta) &\leq  \theta \mu + c \theta^2. 
\end{align}{
Now, for any $\Gamma > 0$ and $\theta > 0$, using standard arguments and \eqref{eqd:1}, we obtain} 
\begin{align}
\Pb(\sum_{i=1}^n X_i \leq n \mu - \Gamma) & = \Pb\big(\exp(-\theta(\sum_{i=1}^n X_i - n \mu)) \geq \exp(\Gamma \theta)\big) \nonumber \\
& \leq \exp(-\theta \Gamma + \theta n \mu) \Eb[\exp(- \theta X)]^n \nonumber \\
& = \exp\big(-\theta \Gamma + \theta n \mu + n M_X(-\theta)\big) \nonumber \\
& \leq \exp\big(-\theta \Gamma + c n \theta^2  \big)
\end{align}
{For any $0 < \theta \leq \Gamma/(2nc)$, }
\begin{align*}
\Pb(\sum_{i=1}^n X_i \leq n \mu - \Gamma) & \leq \exp\big(- \frac{1}{2}\Gamma \theta\big).
\end{align*}
{Using $\Gamma =  n \mu \varepsilon$ and $\theta^* = \min(\delta/2,  \mu \varepsilon/(2c))$, we have}
\begin{align}
\Pb(\sum_{i=1}^n X_i \leq n \mu (1-\varepsilon)) & \leq \exp\big(- \frac{1}{2}n \mu \varepsilon \theta^*\big) % \nonumber \\
%& 
= \exp\big(- C_2(\varepsilon, \delta) n \mu\big), 
\end{align}
{where $C_2(\varepsilon, \delta) = \frac{1}{2}\varepsilon \min(\delta/2,  \mu \varepsilon/(2c))$. 
This completes the proof of Proposition \ref{prop:std}.}

\medskip
\noindent
{\bf Proof of Proposition \ref{prop:error}.} {Given $1 \leq s < r$, $\gamma > 0$ and $\theta > 0$, 
using standard arguments (with the notation that the random variable $X$ has an identical distribution as  $X_i, 1\leq i\leq r+s$)}
\begin{align}
\Pb(\sum_{i=1}^r X_i - \sum_{j=1}^s X_{r+j} \leq -\gamma) & = \Pb(-\theta(\sum_{i=1}^r X_i - \sum_{j=1}^s X_{r+j}) \geq \gamma \theta) \nonumber\\
& \leq \exp(-\theta \gamma) \Eb[\exp(-\theta X)]^r \Eb[\exp(\theta X)]^s. \nonumber
\end{align}
{Using notation and arguments similar to that in the proof of Proposition \ref{prop:std}, we conclude that the above inequality can be bounded 
above, for some $c > 0$ and $\theta \in (-\delta_1, \delta_1)$ for $0 < \delta_1 < \delta$ as}
\begin{align}
\Pb(\sum_{i=1}^r X_i - \sum_{j=1}^s X_{r+j} \leq -\gamma) & \leq \exp(-\theta \gamma + (s-r) \theta \mu + (r+s) c \theta^2).
\end{align}
{For $\theta^* = \min(\frac{\gamma + (r-s)\mu}{2 (r+s) c}, \delta_1/2)$, we obtain }
\begin{align}
\Pb(\sum_{i=1}^r X_i - \sum_{j=1}^s X_{r+j} \leq -\gamma) & \leq  \exp\big(- \frac{1}{2}\theta^* (\gamma + (r-s) \mu)\big).
\end{align}

\subsubsection{Proof of Theorem \ref{thm:randTreeError}.}

{Theorem \ref{thm:randTreeError} assumes that the spreading times have an exponential distribution with (unknown) parameter $\lambda > 0$ for all edges. The underlying graph is a generic random tree, just like that in Theorem \ref{thm:randtree}. We shall crucially utilize the `memory-less' property of the exponential distribution to obtain the exponential error bound on $\limsup_{t\to\infty} \Pb(C^k_t)$ claimed in Theorem \ref{thm:randTreeError}.} 

{To that end, continuing with notations from the proof of Theorem \ref{thm:randtree}, let $T_k = \inf\{t > 0: T(t) = k\}$. By definition, } 
\begin{align}
\limsup_{t\to\infty} \Pb(C^k_t | T_k =\infty) & = 0. 
\end{align}
{Therefore,} 
\begin{align}
\limsup_{t\to\infty} \Pb(C^k_t) & \leq \limsup_{t\to\infty} \Pb(C^k_t | T_k <\infty). 
\end{align}
{Therefore, let us assume that $T_k < \infty$ and we will be interested in $t > T_k$. We shall re-define the index for time as $t' = t - T_k$. 
When $t' = 0$, we have exactly $k$ nodes infected and let them be $v_1,\dots, v_k$, chronologically infected in that order. 
Let $d = \eta_k + 1$ denote the total number of neighbors of $v_k$, let $w_1$ denote the neighbor of $v_k$ on the path connecting $v_k$ and $v_1$ and let
$w_2,\dots, w_d$ be the other neighbors of $v_k$.  Let ${\cal T}^k_1(t')$ (with $T^k_1(t') = | {\cal T}^k_1(t') |$) be the sub-tree rooted at $w_1$ including $v_1$ (and not including $v_k$). Similarly, let ${\cal T}^k_j(t')$ (with $T^k_j(t') = | {\cal T}^k_j(t') |$) be the sub-tree rooted at $w_j$, not including $v_k$ 
for $2\leq j\leq d$. By definition $T^k_1(0) = k-1$ and $T^k_j(0) = 0, ~2\leq j\leq d$.} 

{Let $Z(t')$ be the size of the rumor boundary of the graph at $t'$ and let $\zeta_k = Z(0)$ be the size of the rumor boundary immediately after the $k^{th}$ node is  infected. By definition, $\zeta_k \geq d-1$ as $w_2,\dots, w_d$ are part of the rumor boundary when $t' = 0$. Let $X_1(t'), \dots, X_{\zeta_k}(t')$ be the size of the sub-trees at time $t' = t - T_k$ (for $t \geq T_k$) growing from these $\zeta_k$ rumor boundary nodes. Due to the memory-less property of 
exponential spreading time distribution, it can be argued that $X_1(t'),\dots, X_{\zeta_k}(t')$ are independent and identically distributed random variables. Putting the above discussion together, we have that }
\begin{align}
T_1^k(t') & = k -1 + \sum_{j=1}^{\zeta_k - d + 1} X_j(t'), \nonumber \\
\sum_{i=2}^d T_i^k(t') & = \sum_{j=\zeta_k - d + 2}^{\zeta_k} X_j(t').
\end{align}

%
%
%To obtain the upper bound in Theorem \ref{thm:randTreeError} for $\lim_{t\rightarrow\infty}\Pb\paranth{C^k_{t}}$ we assume that after time $t$ at least $k$ nodes have been infected ($n(t)\geq k$), with the $k^{th}$ infected node being defined as $v_k$ which has degree $d$.  The number of children of $v_k$ is $\eta_k=d-1 \geq 2$.  $C^k_{t}$ is the event that $v_k$ is the rumor center after time $t$.  To upper bound the probability of this event, we will use the memoryless property of the exponential spreading times crucially.  There are $d$ subtrees neighboring $v_k$.  We define the time when $v_k$ is infected as $t_k<t$.  We define the size of the rumor boundary at $t_k$ as $Z(t_k)$ which consists of all uninfected nodes neighboring infected nodes.  We have that $Z(t_k)=1+\sum_{i=1}^{k}\paranth{\eta_i-1}$.  Because of the memoryless property of the exponential spreading times and the way in which the random tree is constructed, at time $t_k$ each node in the rumor boundary is an independent copy of identically distributed subtree random processes which we will refer to as $X_j(\cdot)$, for $1 \leq j\leq Z(t_k)$ in the rumor boundary.  This is illustrated
%in Figure \ref{fig:XrandTreeError}.   
%
%%%%%%%%%%%%%%%%%%%%%%%%%%%%%%%%%%%%%%%%%%%%%%%%%%%%%%%%
\begin{figure}\centering
	\includegraphics[scale=1]{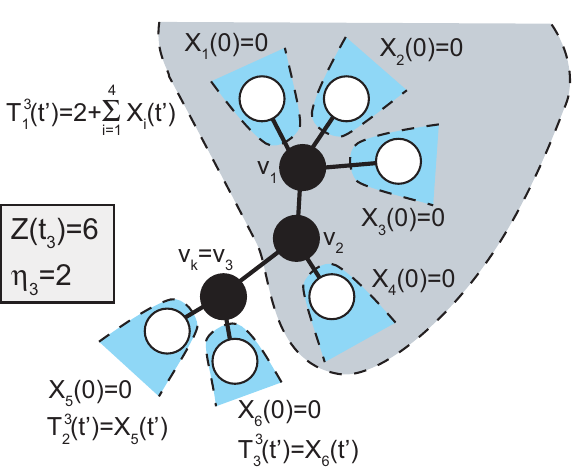}
	\caption{Illustration of the labeling of the
	subtree random processes $X_j(t')$ 
	for $k=3$ in a rumor graph at $t=t_k$ 
	(the time of infection of $v_k$).  
	The rumor infected nodes are colored black, 
	and the uninfected nodes are white.}
	\label{fig:XrandTreeError}
\end{figure}
%%%%
%
%We now use the notation from Section \ref{sssec:Cnk} for the subtree processes.  Specifically, let us imagine the node $v_k$ as 
%the global root and with respect to it, let $T^k_1(t)$ denoted the size of the subtree rooted at $v_{k-1}$ at
%time $t$. Again we define $t'=t-t_k$.  Let  $T^k_i(t')$ 
%for $i=2,...,d$ be the other subtrees rooted at the 
%children of $v_k$ (which were not infected at time 
%$t_k$ but were only part of the rumor boundary). 
%Then, we have that $T^k_1(0)=k-1$ and 
%$T^k_i(0)=0$ for $2\leq i\leq d$.  In $T^k_1(0)$, 
%there are $Z(t_k)-(d-1)$ nodes on the rumor 
%boundary, each of which will be the source for i.i.d. 
%subtree process $X_{j}(t')$, $1\leq j\leq Z(t_k)-(d-1)$ starting at $t'=0$.  
%Therefore,
%\begin{align}
%	T^k_1(t')& = k-1+ \sum_{j=1}^{Z(t_k)-(d-1)}X_{j}(t').
%\end{align}
%and 
%\begin{align}
%	T^k_j(t')& = X_{Z(t_k)-d+j}(t'), ~~~~~2\leq j \leq d.
%\end{align}
%%
%Now, we will upper bound $\Pb\paranth{C^k_{t}}$ \textbf{as follows}.   
{With $T_k < \infty$, for $t \geq T_k$ and $t' = t - T_k$, }
\begin{align}
	\Pb\paranth{C^k_{t} | T_k < \infty}&\leq	\Pb\paranth{\cap_{i=1}^{d} \Big\{2T^k_i(t') \leq \sum_{j=1}^d T^k_j(t')\Big\}+1 \big| T_k < \infty }\nonumber\\
	&\leq \Pb\paranth{T^k_1(t') \leq \sum_{j=2}^{d} T^k_j(t')+1 | T_k < \infty}\nonumber\\
	&\leq \Pb\paranth{k-1+ \sum_{j=1}^{\zeta_k-d+1} X_{j}(t') \leq \sum_{j'=\zeta_k-d+2}^{\zeta_k} X_{j'}(t')+1}
	\label{eq:Pcnk}.
	\end{align}
	%%%
{We shall argue that the term on the right in \eqref{eq:Pcnk} is bounded from above by $O(\exp(- a k))$ for all $k$ large enough. To that end, we shall
utilize Propositions \ref{prop:std} and \ref{prop:error}. }

{First, recall that $\zeta_k - d + 1$, which is the total number of nodes in the rumor boundary at time $t = T_k$ due to the first $k-1$ infected
nodes, equals $\sum_{i=1}^{k-1} (\eta_i - 1)$, where $\eta_1,\dots, \eta_{k-1}$ are the random number of children of the first $k-1$ infected nodes.  By assumption, $\Eb[\eta] > 1$. Therefore, using Proposition \ref{prop:std}, it follows that for an appropriate choice of constants $C_1, C_2$, }
\begin{align}\label{eqd:2}
\Pb\big(\zeta_k - d + 1 \leq (k-1)(\Eb[\eta]-1)/2\big) & \leq C_1 \exp(- C_2 k). 
\end{align} 
{Second, consider the rumor boundary induced due to the children of $v_k$, denoted in the above sum as $d-1$ nodes 
(corresponding to the terms in the right hand side of the equation). Since $d$, the degree of $v_k$ is a random number
distributed as per $\eta$ and $\Eb[\exp(\theta\eta)] < \infty$ for all $\theta \in (-\varepsilon, \varepsilon)$ for some 
$\varepsilon > 0$, it follows that for appropriate constants $C_3, C_4 > 0$ (with $k\geq 2$), }
\begin{align}\label{eqd:3}
\Pb(d-1 > (k-1)(\Eb[\eta]-1)/4) & \leq C_3 \exp(-C_4 k). 
\end{align}
{Define the event $E = \{\zeta_k - d + 1 > (k-1)(\Eb[\eta]-1)/2\} \cap \{d-1 \leq (k-1)(\Eb[\eta]-1)/4\}$. Then, from \eqref{eqd:2}-\eqref{eqd:3}, 
we have $\Pb(E^c) \leq C_5 \exp(- C_6 k)$ where $C_5 = C_1 + C_3$ and $C_6 = \min(C_2, C_4)$.} 

{Finally, to bound $\Pb(F)$, where $F = \{k-1+ \sum_{j=1}^{\zeta_k-d+1} X_{j}(t') \leq \sum_{j'=\zeta_k-d+2}^{\zeta_k} X_{j'}(t')+1\}$, consider
the following: for all $k$ large enough, using Proposition \ref{prop:error}, we have}
\begin{align}
\Pb(F) & \leq \Pb(F | E) + \Pb(E^c) \nonumber \\
 & \leq  C_7 \exp( - C_8 k) +  C_5 \exp(- C_6 k) ~=~ C' \exp(-C'' k).
\end{align}
{In the last inequality, the first term is derived by applying Proposition \ref{prop:error} where $r \geq (k-1)(\Eb[\eta]-1)/2$ and $s \leq (k-1)(\Eb[\eta]-1)/4$, i.e. $r \geq 2 s$, $\gamma = k-2$, and $C', C'' > 0$ are appropriate constants depending upon $C_5, C_6, C_7$ and $C_8$. Note that the conditions of Proposition \ref{prop:error} are satisfied 
because of Theorem \ref{thm:mgf}. This completes the proof of 
Theorem \ref{thm:randTreeError}.}

%%%%%%%%%%%%%%%%%%%%%%%%%%%%%%%%%%%%%%%%%%%%%%%%%%%%%%%%%%%%%%%%%%%%%%%%%%%%%%%%
%%%%%%%%%%%%%%%%%%%%%%%%%%%%%%%%%%%%%%%%%%%%%%%%%%%%%%%%%%%%%%%%%%%%%%%%%%%%%%%%
%%%%%%%%%%%%%%%%%%%%%%%%%%%%%%%%%%%%%%%%%%%%%%%%%%%%%%%%%%%%%%%%%%%%%%%%%%%%%%%%
%%%%%%%%%%%%%%%%%%%%%%%%%%%%%%%%%%%%%%%%%%%%%%%%%%%%%%%%%%%%%%%%%%%%%%%%%%%%%%%%
\subsection{Proof of Theorem \ref{thm:geom}: Geometric Trees}

The proof of Theorem \ref{thm:geom} uses
the characterization of the rumor center provided by 
Proposition \ref{lem:rc}. That is, we wish to
show that for all $n$ large enough, the event that the size of the $d^*$ rumor infected sub-trees 
of the source $v^*$ are 
essentially `balanced' occurs
with high probability. 
To establish this, we shall use coarse estimations
on the size of each of these sub-trees using the standard
concentration property of renewal processes along with
geometric growth. This will be unlike the proof for 
regular trees where we had to necessarily delve into 
very fine detailed probabilistic estimates of the 
size of the sub-trees to establish the result. This 
relatively easier proof for geometric trees (despite their
heterogeneity) brings out the fact that it is 
fundamentally much more difficult to analyze expanding
trees than geometric structures as expanding trees do 
not yield to generic concentration based estimations
as they necessarily have very high variances. 

To that end, we shall start by obtaining sharp estimates on
the size of each of the $d^*$ rumor infected  sub-trees of
$v^*$ for any given time $t$.  We are assuming here that the spreading times 
have a distribution $F$ with mean $\mu>0$ and an exponential
tail (precisely, if $X$ is random variable with distribution $F$, then 
$\Eb[\exp(\theta X)] < \infty$ for $\theta \in (-\beps,\beps)$ for some $\beps > 0$).
 Initially, at time $0$
the source node $v^*$ is infected with the rumor. It starts spreading to 
its $d^*$ children (neighbors). Let $T_i(t)$ denote the size
of the rumor infected subtree, denoted by $\mathcal T_i(t)$, 
rooted at the $i$th child (or neighbor) of node $v^*$. 
Initially, $T_i(0) = 0$.  Due to the balanced and geometric growth 
conditions assumed in Theorem \ref{thm:geom}, the following
will be satisfied: for small enough $\epsilon > 0$
(a) every node within a distance $\frac{t}{\mu}\left(1-\epsilon\right)$ 
of $v^*$ is in one of the $\mathcal T_i(t)$, and (b) no node beyond distance 
$\frac{t}{\mu}\left(1+\epsilon\right)$ of $v^*$ is in any of the $\mathcal T_i(t)$. 
Such a tight characterization of the `shape' of $\mathcal T_i(t)$ along with
the polynomial growth will provide sharp enough bound on $T_i(t)$
that will result in establishing Theorem \ref{thm:geom}. 
This result is summarized below with its proof in Section \ref{sec:fill}.
%
%\begin{figure}
	%\centering
	%	\includegraphics[height=1.5in,width=2in]{figures/fill.eps}
	%	\caption{The structure of a rumor graph on a geometric tree after the rumor has spread for time $t$.}  
	%\label{fig:fill}
%\end{figure}
%%%%%%%%%%%%%%%%%%%%%%%%%%%%%%%%%%%%%%%%%%%%%%%%%%%%%%%%%%%%%%%%%%%%%%%%%%%%%%%%%%%%%%%%%%%%%%%%%%%%%%%%%%%
\begin{prop}\label{prop:fill}
Consider a geometric tree with parameters $\alpha > 0$ and $0< b \leq c$ 
as assumed in Theorem \ref{thm:geom} and let the rumor spread from 
source $v^*$ starting at time $0$ as per the SI model with spreading time
distribution $F$ such that the mean is $\mu$ and 
$\Eb[\exp(\theta X)] < \infty$ for $\theta \in (-\beps, \beps)$ for some $\beps > 0$ where 
$X$ is distributed as per $F$.  Define $\epsilon = t^{-1/2+\delta}$ 
for any  $0< \delta < 1/2$. Let $G(t)$ be the rumor
infected tree at time $t$. Let ${\mathcal G}_t$ be the
set of all trees rooted at $v^*$ (rumor graphs) such that all nodes
within distance $\frac{t}{\mu}(1-\epsilon)$ from $v^*$ are in the tree and no node beyond distance $\frac{t}{\mu}(1+\epsilon)$ from $v^*$ is in the tree.  
Then%
\begin{align*}
\mathbf P(G_t \in \mathcal G_t) & = 1 - O\big(e^{-t^\delta}\big)  ~\stackrel{t\to\infty}{\longrightarrow}~ 1.
\end{align*}
\end{prop}
%%%%%%%%%%%%%%%%%%%%%%%%%%%%%%%%%
Define $\mathcal E_t$ as the event that $G_t \in \mathcal G_t$. Under
event $\mathcal E_t$, consider the sizes of the sub-trees $T_i(t)$ for $1\leq i\leq d_{v^*}$.
Due to the polynomial growth condition and $\mathcal  E_t$, we obtain the following
bounds on each $T_i(t)$ for all $1\leq i\leq d_{v^*}$:
\begin{align*}
\sum_{r=1}^{\frac{t}{\mu}(1-\epsilon)-1} b r^\alpha & \leq T_i(t)% \\
%                                      & 
~\leq \sum_{r=1}^{\frac{t}{\mu}(1+\epsilon)-1} c r^\alpha.
\end{align*}
Now bounding the summations by Riemann integrals, we have
\begin{align*}
\int_0^{L-1} r^\alpha dr & \leq \sum_{r=1}^L r^\alpha ~\leq~\int_0^{L+1} r^\alpha dr. 
\end{align*}
%%%%%%%%%%%%%%%%%%%%%%%%%%%%%%%%%%%%%5
Therefore, it follows that under event $\mathcal E_t$, for all $1\leq i\leq d_{v^*}$
\begin{align*}
\frac{b}{1+\alpha} \paranth{\frac{t}{\mu}(1-\epsilon)-2}^{\alpha+1} & \leq T_i(t) ~\leq~ \frac{c}{1+\alpha}  \paranth{\frac{t}{\mu}(1+\epsilon)}^{\alpha+1}. 
\end{align*}

In the most `unbalanced' situation, $d_{v^*}-1$ of these sub-trees 
have minimal size $T_{\text{min}}(t)$ and the remaining one
sub-tree has size $T_{\text{max}}(t)$ where 

\begin{align*}
T_{\text{min}}(t) & = \frac{b}{1+\alpha} \paranth{\frac{t}{\mu}(1-\epsilon)-2}^{\alpha+1}, \\
T_{\text{max}}(t) & = \frac{c}{1+\alpha}  \paranth{\frac{t}{\mu}(1+\epsilon)}^{\alpha+1}.
\end{align*}

Since by assumption $c < b (d_{v^*}-1)$, there exists $\gamma > 0$ such that $(1+\gamma)c <  b(d_{v^*}-1)$. Therefore, for any choice of $\epsilon = t^{-1/2 + \delta}$ for
some $\delta \in (0,1/2)$, we have

\begin{align*}
\frac{(d^*-1)T_{\text{min}}(t)+1}{T_{\text{max}}(t)} =& \frac{b (d_{v^*}-1)}{c} \paranth{\frac{\frac{t}{\mu} - 											t^{\frac{1}{2}	+\delta} - 2}{\frac{t}{\mu} + t^{\frac{1}{2} + \delta}}}^{\alpha+1} \\
		& + \frac{1+\alpha}{c}\paranth{\frac{1}{\frac{t}{\mu}+ t^{\frac{1}{2} + \delta}}}^{\alpha+1}\\
\stackrel{(i)}{>}& (1+\gamma) \paranth{\frac{1-t^{-\frac{1}{2}+ \delta}\mu - 2\mu t^{-1}} 
				{1 +t^{-\frac{1}{2}+\delta}\mu}}^{\alpha+1} \\
				&+\frac{1+\alpha}{c}\paranth{\frac{1}{\frac{t}{\mu} + t^{\frac{1}{2} + \delta}}}^{\alpha+1}\\
 > & 1+\gamma\\
 >& 1,
\end{align*}
for $t$ large enough since as $t\to\infty$ the first term in inequality (i) 
goes to $1$ and the second term goes to $0$.  From this, it immediately follows
that under event $\mathcal E_t$ for $t$ large enough

\[ \max_{1\leq i\leq d_{v^*}} T_i(t) < \frac{1}{2} \left(\sum_{i=1}^{d_{v^*}} T_i(t)+1\right).\]

Therefore, by Lemma \ref{lem:rc} it follows that the rumor
center is unique and equals $v^*$.  We also have that $\mathcal E_t\subset C^1_{t}$.  Thus, from above and Proposition \ref{prop:fill} we obtain

\begin{align*}
\lim_{t} \Pb\big(C^1_{t}\big) & \geq \lim_{t} \Pb\big({\mathcal E}_t\big) \\
                                       & = 1.
\end{align*}

%%%%%%%%%%%%
This completes the proof of Theorem \ref{thm:geom}.

%%%%%%%%%%%%%%%%%%%%%%%%%%%%%%%%%%%%%%%%%%%%%%%%%%%%%%%%%%%%%%%%%%%%%%%%%%%%%%
%%%%%%%%%%%%%%%%%%%%%%%%%%%%%%%%%%%%%%%%%%%%%%%%%%%%%%%%%%%%%%%%%%%%%%%%%%%%%%
%%%%%%%%%%%%%%%%%%%%%%%%%%%%%%%%%%%%%%%%%%%%%%%%%%%%%%%%%%%%%%%%%%%%%%%%%%%%%%
%\section*{Acknowledgements}
%\appendix

\subsubsection{Proof of Proposition \ref{prop:fill}}\label{sec:fill}
We recall that Proposition \ref{prop:fill} stated that for a rumor spreading for time $t$ as per the SI model with a general distribution with mean spreading time $\mu$ the rumor graph on a geometric tree is full up to a distance $\frac{t}{\mu}(1-\epsilon)$ from the source and does not extend beyond $\frac{t}{\mu}(1+\epsilon)$, for $\epsilon = t^{-1/2+\delta}$ for some positive 
$\delta \in (0,1/2)$. To establish this, we shall use the following well 
known concentration property of renewal processes. We 
provide its proof later for completeness. 
%%%%%%%%%%%%%%%%%%%%%%%%%%%%%%%%%%%%%%%%%%%%%%%%%%%%%%%%%%%%%%%%%%%%%%%%%%%%%%%%%%%%%%%%%%%%%%%%%%%%%%%%%%%%%%%%%%%
\begin{prop}\label{prop:poisson}
Consider a renewal process $P(\cdot)$ with holding times with mean $\mu$ and finite moment generating function in the interval $(-\beps,\beps)$ for some $\beps > 0$. 
Then for any $t > 0$ and any $\gamma \in (0, \beps')$ for a small enough $\beps' > 0$,  there exists a positive constant $c$ such that 
\begin{align*}
	\Pb\paranth{\left|P(t)-\frac{t}{\mu}\right|\geq \frac{t\gamma}{\mu}}\leq 2e^{-\frac{\gamma^2\mu }{8c}t}
\end{align*}
\end{prop}
%%%%%%%%%%%%%%%%%%%%%%%%%%%%%%%%%%%%%%%%%%%%%%%%%%%%%%%%%%%%%%%%%%%%%%%%%%%%%%%%%%%%%%%%%%%%%%%%%%%%%%%%%%%%%%%%%%%
Now we use Proposition \ref{prop:poisson} to establish Proposition \ref{prop:fill}. Recall that
the spreading time along each edge is an independent and identically distributed random variable with mean $\mu$. Now the underlying network graph is a tree. Therefore
for any node $v$ at distance $r$ from source node $v^*$, there is a unique
path (of length $r$) connecting $v$ and $v^*$. Then, the spread of the rumor
along this path can be thought of as a renewal process, say $P(t)$, 
and node $v$ is infected by time $t$ if and only if $P(t) \geq r$. Therefore,
from Proposition \ref{prop:poisson} it follows that for any node $v$ that is
at distance $\frac{t}{\mu}(1-\epsilon)$ for $\epsilon = t^{-\frac{1}{2} + \delta}$ for some
$\delta \in (0,1/2)$ (for all $t$ large enough), 
%%%%%%%%%%%%%%5
\begin{align*}
\Pb\big(v \text{~is not rumor infected}\big) & \leq 2 e^{-\frac{\epsilon^2\mu t}{8c}} \\
                                            & = 2 e^{-\frac{\mu}{8c}t^{2\delta}}. 
\end{align*}
Now the number of such nodes at distance $\frac{t}{\mu}(1-\epsilon)$ from $v^*$ is 
at most $O\paranth{\paranth{\frac{t}{\mu}}^{\alpha+1}}$ (which follows from arguments similar to those in
the proof of Theorem \ref{thm:geom}). Therefore, by an application of the union
bound it follows that 
\begin{align*}
& \Pb\paranth{\text{a node at distance $\frac{t}{\mu}(1-\epsilon)$ from $v^*$ isn't infected}} \\
& \qquad = O\paranth{2\paranth{\frac{t}{\mu}}^{\alpha+1} e^{-\frac{\mu}{8c}t^{2\delta}}} \\
& \qquad = O\paranth{e^{-\frac{\mu}{8c}t^{\delta}}}.
\end{align*}
%for $t$ large enough. 
Using similar argument and another 
application of Proposition \ref{prop:poisson}, it can be argued that
\begin{align*}
& \Pb\paranth{\text{a node at distance $t(1+\epsilon)$ from $v^*$ is infected}} \\
& \qquad = O\paranth{e^{-\frac{\mu}{8c}t^{\delta}}}. 
\end{align*}
Since the rumor is a `spreading' process, if all nodes at distance $r$ from $v^*$ are
infected, then so are all nodes at distance $r' < r$ from $v^*$; if all nodes
at distance $r$ from $v^*$ are not infected then so are all nodes at distance $r' > r$ 
from $v^*$. Therefore, it follows that with probability $1-O\paranth{e^{-\frac{\mu}{8c}t^{\delta}}}$, 
all nodes at distance up to $\frac{t}{\mu}(1-\epsilon)$ from $v^*$ are infected and all
nodes beyond distance $\frac{t}{\mu}(1+\epsilon)$ from $v^*$ are not infected. This
completes the proof of Proposition \ref{prop:fill}.

%%%%%%%%%%%%%%%%%%%%%%%%%%%%%%%%%%%%%%%%%%%%%%%%%%%%%%%%%%%%%%%%%%%%
\subsubsection{Proof of Proposition \ref{prop:poisson}} \label{sec:poisson}

We wish to provide bounds on the probability of $P(t) \leq \mu t(1-\gamma)$ and
$P(t) \geq \mu t (1+\gamma)$ for a renewal process $P(\cdot)$ with holding times with mean $\mu$ and finite moment generating function. Define the $n^{th}$ arrival time $S_n$ as
\begin{align*}
	S_n & = \sum_{i=1}^{n}X_i
\end{align*}
where $X_i$ are non-negative i.i.d. random variables with a well defined moment generating function $M_X(\theta) = \Eb[\exp(\theta X)] < \infty$ for $\theta \in (-\beps, \beps)$ for some $\beps > 0$ 
and mean $\Eb\bracket{X_i}=\mu>0$.  We can relate the arrival times to the renewal process by the following relations:
\begin{align*}
	\Pb\paranth{P(t)\leq n}&= \Pb\paranth{S_n\geq t}
\end{align*}
and
\begin{align*}
	\Pb\paranth{P(t)\geq n}&= \Pb\paranth{S_n\leq t}.
\end{align*}
The first relation says that the probability of less than $n$ arrivals in time $t$ is equal to the probability that the $n$th arrival happens after time $t$.  The second relation says that the probability of more than $n$ arrivals in time $t$ is equal to the probability that the $n$th arrival happens before time $t$.

We now bound 	$\Pb\paranth{S_n\geq t}$.  To that end, for $\theta \in (0,\beps)$ it follows from the Chernoff bound that 
\begin{align*}
	\Pb\paranth{S_n\geq t}& = 	\Pb\paranth{e^{\theta S_n}\geq e^{\theta t}}\\
	&\leq M_X\paranth{\theta}^n e^{-\theta t}.
\end{align*}
%%%%%%%%%%%%%%%%%%%%%%%%%%%%%%%%%%%%%%%%%%%%%%%%%%%%%%%%%%%%%%%%%%%%%%%%%%%%
We can use the following approximation for $M_X\paranth{\theta}$ which is valid for small $\theta$, say $\theta \in (0,\beps^+)$ for $0 < \beps^+ \leq \beps$.
\begin{align*}
	M_X\paranth{\theta}&= 1+\theta\mu+ \theta^2\frac{\Eb\bracket{X^2}}{2}+ \theta^3\sum_{i=3}^{\infty}\theta^{i-3}\frac{\Eb\bracket{X^i}}{i!}\\
	&\leq 1+\theta\mu+ c_1\theta^2
\end{align*}
for some finite positive constant $c_1$.  Using this along with the inequality $\log\paranth{1+x}\leq x$ for $-1<x$, we obtain
\begin{align*}
	\log\paranth{\Pb\paranth{S_n\geq t}}&\leq                                                                  n\log\paranth{M_X\paranth{\theta}}-\theta t\\
	&\leq n\log\paranth{1+\theta\mu+ c_1\theta^2}-\theta t\\
	&\leq \theta\paranth{\mu n-t}+nc_1\theta^2.
\end{align*}
%%%%%%%%%%%%%%%%%%%%%%%%%%%%%%%%%%
To minimize this probability, we find the $\theta$ that minimizes $\theta\paranth{\mu n-t}+nc_1\theta^2$.  This happens for $\theta = \frac{1}{2c_1}\paranth{\frac{t}{n}-\mu}$.  We set $n = \frac{t}{\mu}\paranth{1-\gamma}$, so the minimum value is achieved for $\theta^* = \frac{\gamma\mu}{2c_1(1-\gamma)}$.  Therefore, there exists $\beps_1 > 0$ so that for $\gamma \in (0,\beps_1)$, the corresponding $\theta^* = \frac{\gamma\mu}{2c_1(1-\gamma)} < \beps^+$, so that
the quadratic approximation of $M_X(\theta)$ is valid. Given this,  we obtain
%%%%%%%%%%%%%%%%%%%%%%%%%%%%%%%%%%%%%%%%%%%%%%%%%%%%%%%%
\begin{align*}
	\log\paranth{\Pb\paranth{S_{\frac{t}{\mu}(1-\gamma)}\geq t}} &\leq -\frac{\gamma\mu}{2c_1(1-\gamma)}\paranth{\gamma t}+\frac{tc_1}{\mu}(1-\gamma) \frac{\gamma^2\mu^2}{4c_1^2\paranth{1-\gamma}^2}\\
	&\leq -\frac{\gamma^2\mu t}{2c_1(1-\gamma)}+\frac{\gamma^2\mu t}{4c_1(1-\gamma)}\\
	&\leq -\frac{\gamma^2\mu t}{4c_1(1-\gamma)}\\
	&\leq -\frac{\gamma^2\mu t}{8c_1}.
\end{align*}
%%%%%%%%%%%%%%%%%%%%%%%%%%%%%%%%%%%%%%%%%%%%%%%
With this result, we obtain
\begin{align*}
	\Pb\paranth{P(t)\leq \frac{t}{\mu}\paranth{1-\gamma}}\leq e^{-\frac{\gamma^2\mu t}{8c_1}}, 
\end{align*}
for any $t$ and $\gamma \in (0,\beps_1)$. 
%%%%%%%%%%%%%%%%%%%%%%%%%%%%%%%%%%%%%%%%%%
For the upper bound, we have for $\theta > 0$ 
%%%%%%%%%%%%%%%%%%%%%%%%%%%%%%%%%%%%
\begin{align*}
	\Pb\paranth{S_n\leq t}& = 	\Pb\paranth{e^{-\theta S_n}\geq e^{-\theta t}}\\
	&\leq M_X\paranth{-\theta}^n e^{\theta t}.
\end{align*}
%%%%%%%%%%%%%%%%%%%%%%%%%%%%%%%%%%%%%%%%
We can use the following approximation for $M_X\paranth{-\theta}$ which is valid for small enough $\theta\in (0,\beps^-)$ with $0< \beps^- \leq \beps$.
\begin{align*}
	M_X\paranth{-\theta}&= 1-\theta\mu+ \theta^2\frac{\Eb\bracket{X^2}}{2}- \theta^3\sum_{i=3}^{\infty}\theta^{i-3}\paranth{-1}^{i-3}\frac{\Eb\bracket{X^i}}{i!}\\
	&\leq 1-\theta\mu+ c_2\theta^2
\end{align*}
for some finite positive constant $c_2$.  Using this we obtain
\begin{align*}
	\log\paranth{\Pb\paranth{S_n\leq t}}&\leq                                                                  n\log\paranth{M_X\paranth{-\theta}}+\theta t\\
	&\leq n\log\paranth{1-\theta\mu+ c_2\theta^2}+\theta t\\
	&\leq \theta\paranth{t-\mu n}+nc_2\theta^2.
\end{align*}
%%%%%%%%%%%%%%%%%%%%%%%%%%%%%%%%%%%%%%%%%%%
To minimize this probability, we find the $\theta$ that minimizes $\theta\paranth{t-\mu n}+nc_2\theta^2$.  This happens for $\theta = \frac{1}{2c_2}\paranth{\mu-\frac{t}{n}}$.  We set $n = \frac{t}{\mu}\paranth{1+\gamma}$, so the minimum value is achieved for $\theta^* = \frac{\gamma\mu}{2c_2(1+\gamma)}$. There exists, $\beps_2 > 0$ so that for all $\gamma \in (0,\beps_2)$, $\theta^* = \frac{\gamma\mu}{2c_2(1+\gamma)} \leq \beps^-$ and thus guaranteeing the
validity of quadratic approximation of $M_X(-\theta)$ that we have assumed. Subsequently, we obtain %Using this value, we obtain
%%%%%%%%%%%%%%%%%%%%%%%%%%%%%%%%%%%%%%%%%%%%%%%%%%%%%%%%
\begin{align*}
	\log\paranth{\Pb\paranth{S_{\frac{t}{\mu}(1+\gamma)}\leq t}} &\leq -\frac{\gamma\mu}{2c_2(1+\gamma)}\paranth{\gamma t}+\frac{tc_2}{\mu}(1+\gamma) \frac{\gamma^2\mu^2}{4c_2^2\paranth{1+\gamma}^2}\\
	&\leq -\frac{\gamma^2\mu t}{2c_2(1+\gamma)}+\frac{\gamma^2\mu t}{4c_2(1+\gamma)}\\
	&\leq -\frac{\gamma^2\mu t}{4c_2(1+\gamma)}\\
	&\leq -\frac{\gamma^2\mu t}{8c_2}.
\end{align*}
%%%%%%%%%%%%%%%%%%%%%%%%%%%%%%%%%%%%%%%%%%%%%%%
With this result, we obtain
\begin{align*}
	\Pb\paranth{P(t)\geq \frac{t}{\mu}\paranth{1+\gamma}}\leq e^{-\frac{\gamma^2\mu t}{8c_2}},
\end{align*}
for any $t$ and $\gamma \in (0,\beps_2)$. 

If we set $c = \max\paranth{c_1,c_2}$ and $\beps' = \min(\beps_1,\beps_2)$ and combine the upper and lower bounds then we obtain
\begin{align*}
	\Pb\paranth{\left|P(t)-\frac{t}{\mu}\right|\geq \frac{t\gamma}{\mu}}\leq 2e^{-\frac{\gamma^2\mu }{8c}t}, 
\end{align*}
for any $t$ and $\gamma \in (0,\beps')$ with $\beps' > 0$. This completes the proof of Proposition \ref{prop:poisson}.

%%%%%%%%%%%%%%%%%%%%%%%%%%%%%%%%%%%%%%%%%%%%%%%%%%%%%%%%%%%%%%%%%%%%%%%%%%%%%%%%%%%
\section{Conclusion}

Finding the source of a rumor in a network is an important and challenging problem in many different fields.  Here we characterized the performance of the rumor source estimator known as rumor centrality for generic tree graphs. Our analysis was based upon continuous time branching processes and generalized Polya's urn models.  As an implication of this novel
analysis method, we recovered all the previous results for regular trees from \cite{ref:rc} as a special case.  We also showed that for rumor spreading on a random regular graphs, the probability that the estimated source is more than $k$ hops away from the true source decays exponentially in $k$.  Additionally, we showed that for general random trees and hence for sparse random graphs like Erdos-Renyi graphs, there is a strictly positive probability of correct rumor source detection.  Thus, even though rumor centrality is an ML estimator only
for a very specific setting, it is still very effective for a wide range of other graphs and spreading models.
In summary, we have established the \emph{universality} of rumor centrality as a source estimator across a variety of tree structured graphs and spreading time distributions.

\section*{Acknowledgments} Devavrat Shah would like to acknowledge conversations with David Gamarnik and Andrea Montanari at the Banff Research Institute that seeded this line of work. Both 
authors would like to acknowledge the support of the AFOSR Complex Networks Project, the Army Research Office MURI Award W911NF-11-1-0036 on Tomography of Social Networks and the MIT-Shell Graduate Student Fellowship.

\bibliographystyle{plainnat}
\bibliography{Arxivaap}
\end{document}